\documentclass[final]{siamltex}
\usepackage{psfrag}
\usepackage[normal]{subfigure}
\usepackage{epsfig}
\usepackage{latexsym}
\usepackage{amsmath}
\usepackage{amssymb}

\catcode`\@=11
\@addtoreset{equation}{section}
\renewcommand{\theequation}{\thesection.\arabic{equation}}

\newcommand{\beq}{\begin{equation}}
\newcommand{\eeq}{\end{equation}}
\newcommand{\sech}{\mbox{sech}}

\newcommand{\lam}{\lambda}
\newcommand{\hlam}{\hat{\lambda}}

\def\al{\alpha}
\def\be{\beta}
\def\eps{\varepsilon}
\def\ga{\gamma}
\def\Ga{\Gamma}
\def\la{\lambda}
\def\La{\Lambda}
\def\tla{\tilde{\lambda}}
\def\pa{\partial}
\def\si{\sigma}
\def\C{{\cal C}}
\def\D{{\cal D}}
\def\M{{\cal M}}
\def\O{{\cal O}}

\newtheorem{Theorem}{Theorem}

\newtheorem{Lem}[Theorem]{Lemma}

\newtheorem{Rem}[Theorem]{Remark}

\title{Destabilization of fronts in a class of bi-stable systems}
\author{Arjen Doelman\footnotemark[1]\ \footnotemark[2] , David
  Iron{\footnotemark[1]} and Yasumasa Nishiura\footnotemark[3]}

\begin{document}
\maketitle \renewcommand{\thefootnote}{\fnsymbol{footnote}}
\footnotetext[1]{Korteweg-deVries Instituut, Universiteit van
  Amsterdam, Plantage Muidergracht 24, 1018 TV Amsterdam, the
  Netherlands} \footnotetext[2]{Centrum voor Wiskunde en Informatica,
  P.O. Box 94079, 1090 GB Amsterdam, the Netherlands}
\footnotetext[3]{Laboratory of Nonlinear Studies and Computation,
  Research Institure for Electronic Science, Hokkaido University,
  Kita-ku, Sapporo, 060, Japan}
\renewcommand{\thefootnote}{\arabic{footnote}}
\begin{abstract}
In this article,
we consider a class of bi-stable reaction-diffusion equations in two
components on the real line. We assume that the system is
singularly perturbed, i.e. that the ratio of the diffusion
coefficients is (asymptotically) small. This class admits front
solutions that are asymptotically close to the (stable) front
solution of the `trivial' scalar bi-stable limit system
$u_t = u_{xx} + u(1-u^2)$. However, in the system these fronts
can become unstable by varying parameters. This destabilization
is either caused by the essential spectrum associated to the
linearized stability problem, or by an eigenvalue that 
exists near the essential spectrum. We use the Evans function 
to study the various
bifurcation mechanisms and establish an explicit connection between the
character of the destabilization and the possible appearance of
saddle-node bifurcations of heteroclinic orbits in the existence
problem.    
\end{abstract}



\pagestyle{myheadings} \thispagestyle{plain} \markboth{A. DOELMAN, D.
  IRON AND Y. NISHIURA}{DESTABILIZATION OF FRONTS IN A CLASS OF
  BI-STABLE SYSTEMS}

\section{Introduction}
\label{sec:intro}
\setcounter{Theorem}{0}
\setcounter{equation}{0}

The class of bi-stable reaction-diffusion  equations we consider in
this paper is given by
\begin{equation}
\label{most}
\left\{ 
\begin{array}{rcrcl}
U_t & = & \eps^2 U_{xx} & + & (1+V-U^2)U \\
\tau V_t & = &  V_{xx} & + & F(U^2,V; \eps),
\end{array}
\right.  
\end{equation}
where $F(U^2,V; \eps)$ is a smooth function of $U^2$, $V$ and $\eps$ such
that $F(1,0; \eps) \equiv 0$ and $\lim_{\eps \to 0} F(U^2,V;\eps)$
exists ; $\tau > 0$ is a parameter.  Thus, the system is such that the
background state $(U,V) \equiv (\pm 1,0)$ is always a solution.
We furthermore assume that the ratio of the two diffusion
coefficients, $\eps^2$, is asymptotically small, thus the problem has
a singularly perturbed nature.  We consider the system on the
(unbounded) line, i.e.  $(U,V) = (U(x,t),V(x,t))$ with $(x,t) \in
\mathbb{R} \times \mathbb{R}^+$.  Note that (\ref{most}) is (by
construction) symmetric under
\begin{equation}
 \label{symmU}
 U \to - U.  
\end{equation} 
To motivate the structure of (\ref{most}) we introduce the
fast variable
\begin{equation}
\label{fastvar}
\xi = \frac{x}{\eps},
\end{equation}
so that (\ref{most}) can be written in its equivalent `fast'  form
\begin{equation}
\label{most_f}
\left\{ 
\begin{array}{rcrcl}
U_t & = & U_{\xi\xi} & + & (1+V-U^2)U \\
\eps^2 \tau V_t & = &  V_{\xi\xi} & + & \eps^2 F(U^2,V; \eps).
\end{array}
 \right.  
\end{equation}
Since $U(x,t)$ and $V(x,t)$ are a priori supposed to be bounded on the
entire domain $\mathbb{R} \times \mathbb{R}^+$, we find in the natural
(fast reduced) limit, i.e. $\eps \to 0$ in (\ref{most_f}), that $V
\equiv V_0$ and that $U$ is a solution of the well-studied, scalar
(standard) bi-stable or Nagumo equation,
\begin{equation}
\label{scalarlim}
U_t = U_{\xi\xi} + (1+V_0-U^2)U.  
\end{equation} 
In this paper we interpret the original system, (\ref{most}) or
(\ref{most_f}), as a scalar bi-stable Nagumo equation
(\ref{scalarlim}) in which the coefficient of the linear term is
allowed to evolve by reaction and diffusion on a long, or slow,
spatial scale.  Note that the (slow) dynamics of the $V$-component are
allowed to be completely general, except that it is assumed that the
full system conserves the symmetry (\ref{symmU}) and the background
states $U \equiv \pm 1$, at $V \equiv 0$, of the scalar limit (see
also Remark \ref{rem:moregen}).  A priori, one expects that the
$V$-component of front-like solutions will remain small ($\O(\eps)$)
due to the `boundary conditions' $V = 0$ at $\pm \infty$, so that the
effect of the slowly varying $V(x,t)$-component cannot have a
significant influence on the (well-understood) dynamics of the scalar
Nagumo equation. An important motivation of the research in this paper
is to find out whether or not this intuition is correct.

We will focus completely on the existence and stability issues
associated to the persistence of the asymptotically stable stationary
front solutions of the bi-stable equation (\ref{scalarlim}) with $V_0
= 0$. In fact, this paper can also be seen as a first step towards
analyzing the dynamics (and possibly defects) of striped patterns in a
class of relatively simple bi-stable reaction-diffusion equations,
i.e. (\ref{most}) for $(U,V) = (U(x,y,t), V(x,y,t))$ with $(x,y) \in
\mathbb{R}^2$.  The methods and techniques developed in this paper are
supposed to carry over to the analysis of the existence and stability
of spatially periodic solutions of (\ref{most}) and their
two-dimensional counterparts (the planar fronts and the stripe
patterns).  See also section \ref{sec:simdis}.

The problem of the persistence of the stable front solution of the scalar
bi-stable equation (\ref{scalarlim}) is quite subtle, as can be
expected in the light of recent results on the stability of pulses in
singularly perturbed reaction-diffusion equations of the Gray-Scott
and Gierer-Meinhardt type \cite{dgk2,dgk3}. Such systems can also be
written in the form (\ref{most_f}), however, the scalar limit systems
are mono-stable, i.e. in essence of the form 
$U_t = U_{\xi\xi} - U +
U^2$.  The pulses correspond in this (fast reduced) limit to the
stationary homoclinic solution of $u_{\xi\xi} - u + u^2 = 0$. Thus,
one would expect that the pulses of the full system cannot be stable,
since the stability problem associated to the homoclinic solution has
an $\O(1)$ unstable eigenvalue. Nevertheless, stable pulses of this
type do exist in the Gray-Scott and the Gierer-Meinhardt equation
\cite{dgk2,dgk3}. On the other hand, the stability of the pulses in
these mono-stable equations is strongly related to the freedom one has
in these systems to scale the magnitude of the pulses, i.e. the
amplitude of the stable pulses is asymptotically large in $\eps$ in these
mono-stable systems. Such 
scalings are not possible for the fronts in the bi-stable case, since 
the background states $(\pm 1,0)$ are fixed (and $\O(1)$).

In the analysis of the front solutions, we will find that it is
natural to decompose $F(U^2,V;\eps)$ into a component that has a
factor of $(1+V-U^2)$ and a rest term $G(V;\eps)$ that does not depend
on $U^2$. Hence, we write (\ref{most}) as,
\begin{equation}
\label{decomp}
\left\{ 
\begin{array}{rcrcl}
U_t & = & \eps^2 U_{xx} & + & (1+V-U^2)U \\
\tau V_t & = & V_{xx} & + & (1+V-U^2)H(U^2,V; \eps) + G(V;\eps),
\end{array}
\right.
\end{equation}
with $G(0,\eps) \equiv 0$. Note that this decomposition induces no
restriction on $F(U^2,V;\eps)$ since we have assumed that $F$ is
smooth. In fact
\[
G(V) = F(1+V,V) \; \; {\rm and} \; \; (1+V-U^2)H(U^2,V) =
F(U^2,V)-F(1+V,V).
\]
We will find that the quantities $\frac{\pa G}{\pa V}(0;\eps)$ and
$H(1,0;\eps)$ have a crucial impact on the structure and the dynamics
of the front-like solutions. Therefore, we define
\begin{equation}
\label{def:G1H0}
G_1(\eps) = \frac{\pa G}{\pa V}(0;\eps) \; \; {\rm and} \; \; H_0 = H(1,0;\eps);
\end{equation}
$G_1$ is the main bifurcation parameter used in this paper. Throughout
this paper we assume that $H(U^2,V)$ is non-degenerate, i.e. that
$H(1+V,V)$ is not identically $0$, and that $\tau = \O(1)$
(see Remark \ref{rem:degH}).

In section 2 we will show that as long as $G_1 < 0$ and $\O(1)$, the
front solutions of (\ref{scalarlim}) with $V_0 = 0$ persist in a
regular fashion, in the sense that the system (\ref{most}) has front
solution with $U$-components that are asymptotically and uniformly
close to a front in (\ref{scalarlim}) with $V_0 = 0$, and with
$V$-components that are asymptotically and uniformly small (Theorem
\ref{th:ex_reg}).  However, if $G_1$ becomes $O(\eps^2)$ these fronts
become truly singular, in the sense that $V$ becomes $\O(1)$, while
the $U$-component is close to a front of (\ref{scalarlim}) with $V_0
\neq 0$, on the fast spatial scale (and it converges to $U = \pm 1$ on
the slow spatial scale). Moreover, the front solutions are no longer
uniquely determined, there can be several types of heteroclinic front
solutions if $G_1 = \O(\eps^2)$ that may or may not merge in
saddle-node bifurcations of heteroclinic orbits when $G_1$ is varied
(Theorems \ref{th:ex_suslo} and \ref{th:gen_suslo}).  It should
be noted here that we for simplicity consider $G(V) = -\eps^2 \ga V$
in (\ref{most}) in the singular limit $G_1 = \O(\eps^2)$ althrough
this paper -- see Remark \ref{rem:geng1}. We refer to Figure 1 for a
numerical representation of a regular front (Figure 1a) and a singular
front (Figure 1b).
\begin{figure}[ht]
\begin{center}
  \psfrag{x}{$x$} 
  \subfigure[Regular front, $G_1=-1.0$.]{
      \includegraphics[scale=.45]{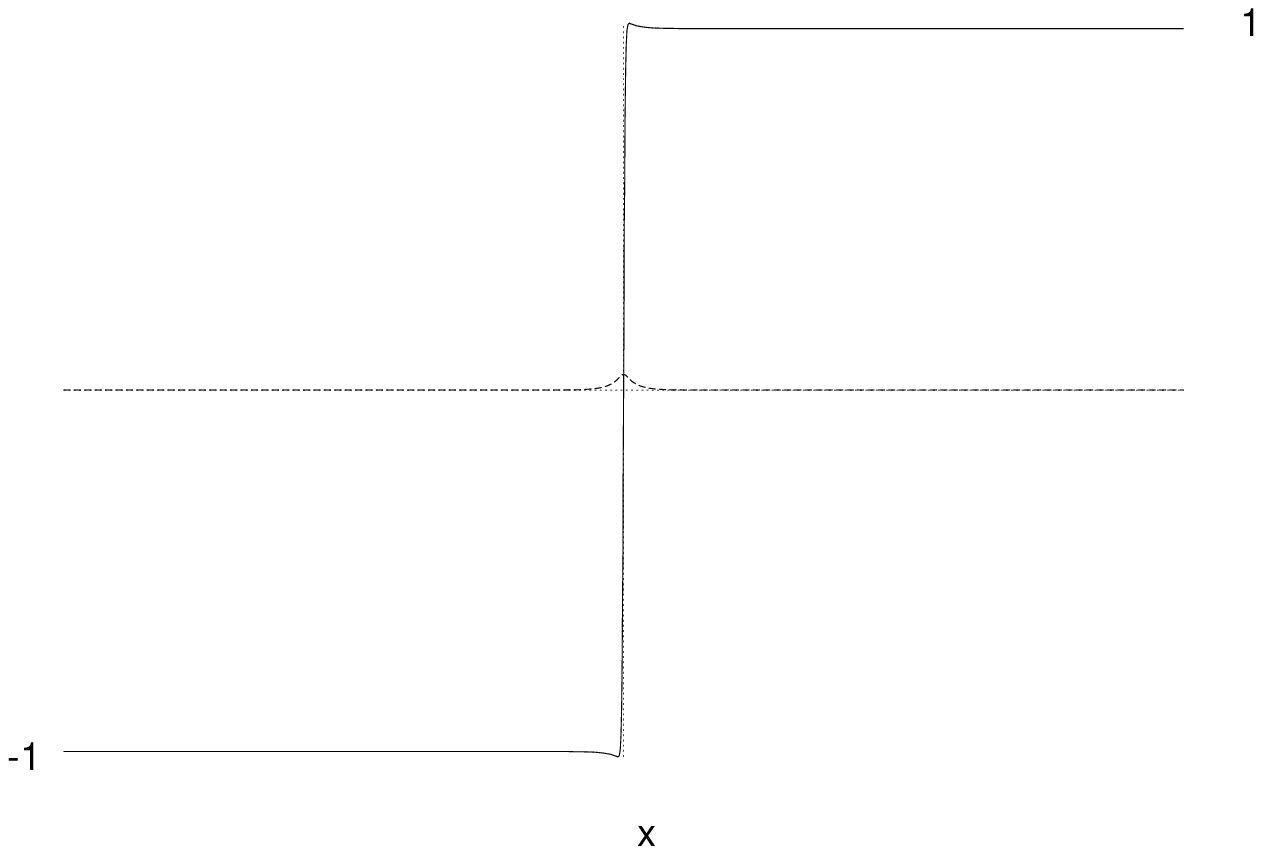}}
  \hspace{2cm} \subfigure[Singular front, $G_1 = -2 \eps^2$.]{
       \includegraphics[scale=.45]{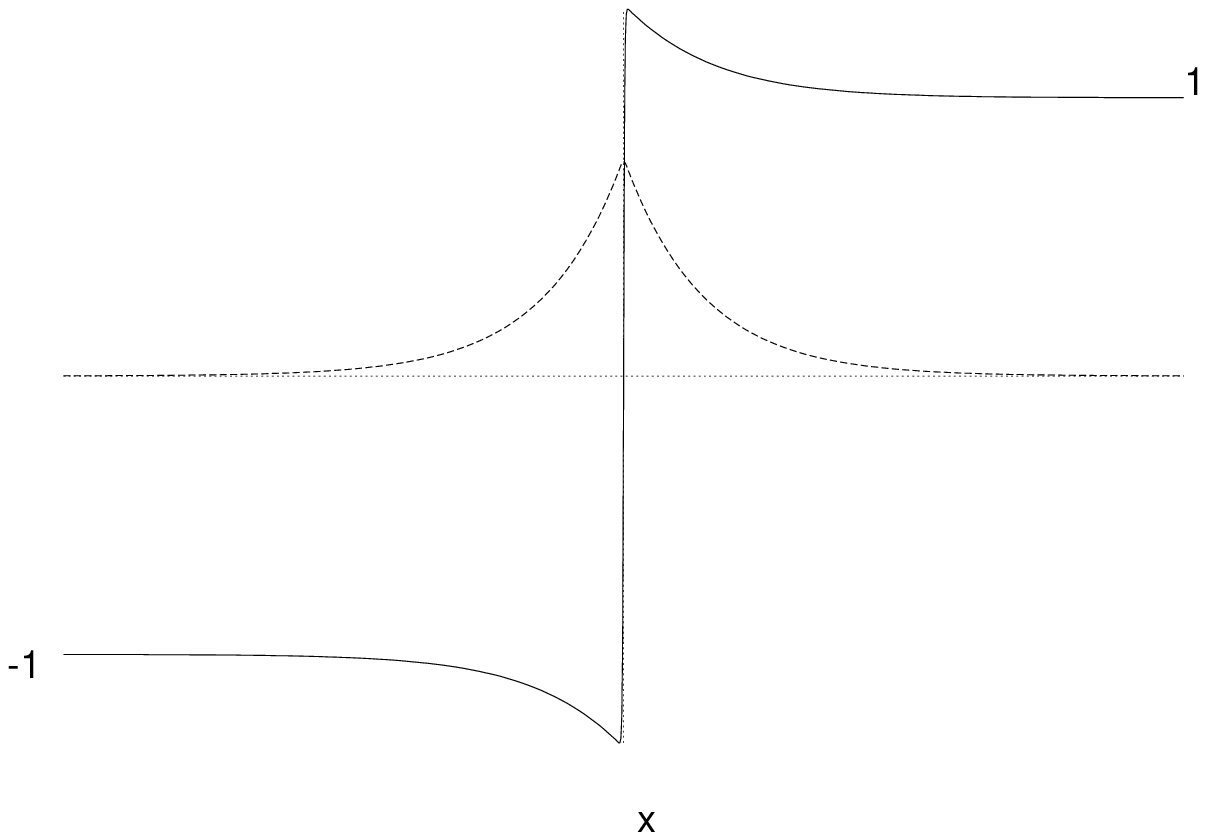}}
\caption{Two asymptoticallly stable front solutions of (\ref{most})/(\ref{decomp})
  plotted on the slow spatial scale $x$ (by a numerical simulation). Here 
  $H(U^2,V) = H_0 U^2$, $G(V) = G_1V$, $\eps=0.1$ and $H_0=1$. 
  The solid curves represent the $U$-coordinates, the dotted
  curves the $V$-coordinates.}
\end{center}
\end{figure}
The magnitude of $G_1$ is also extremely relevant in the stability
analysis.  It can be shown that the (regular) front solutions are
asymptotically stable as long as $G_1 < 0$ and $\O(1)$ and $H_0 + G_1
- 2 \tau < 0$ and $\O(1)$ -- Theorem \ref{th:reg_stab}.  It seems, at
leading order, that the destabilization of the front is caused by the
essential spectrum, $\si_{\rm ess}$, associated to the stability of
the front ($\si_{\rm ess}$ reaches the imaginary axis exactly at $G_1
= 0$ or at $H_0 + G_1 - 2 \tau = 0$ -- Lemma \ref{lem:stabsig_ess}).
However, the analysis also shows that there can be eigenvalues near
the `tips' of $\si_{\rm ess}$, and that it is possible that the
destabilization is caused by such an eigenvalue, i.e. by an element of
the discrete spectrum and not by $\si_{\rm ess}$. These `new'
eigenvalues do not have counterparts in the (scalar) fast reduced
limit problem, they have a singular slow-fast nature and may appear
through edge bifurcations from the essential spectrum.

In section 4 we study in detail the nature of the destabilization as
$G_1<0$ increases towards $0$. In this section it becomes clear that
there is an intimate relation between the geometrical character of the
singularly perturbed existence problem and the character of the
destabilization of the front. This is a natural and frequently
encountered relation -- see for instance \cite{kap98} and the
references there. We establish that a front solution destabilizes at a
critical value of $G_1 = - \eps^2 \ga_{\rm double} < 0$ by an
eigenvalue if and only if it merges with another front solution in a
saddle-node bifurcation of heteroclinic orbits.  Moreover, we are able
to determine the explicit value of this bifurcation value 
$\ga_{\rm double} > 0$.  
If the front does not `encounter' such a saddle-node
as $G_1$ increase to $0$, the front will be destabilized by 
$\si_{\rm ess}$ at $G_1 = 0$ -- see Theorems \ref{th:edge_SN_ex} and
\ref{th:destab}.

Another way to motivate the analysis of this paper is as follows. In
this paper we show that the technique of decomposing the Evans
function associated to the stability of a `localized structure' (a
(traveling) pulse or front) into the product of an analytic `fast' and
a meromorphic`slow' transmission function (\cite{dgk2,dgk3}) can be
extended to a class of bi-stable equations. We show that the slow
transmission function ($t_2(\la, \eps)$) is a natural tool for
analyzing the existence or appearance of eigenvalues near or from the
essential spectrum, and that such eigenvalues play a crucial role in
the stability of the front. Note that in this sense, the theme of this
paper is similar to that of \cite{ks}, where Evans function techniques
are developed to study eigenvalues near $\si_{\rm ess}$ in a class of
nearly integrable systems.

The paper is organized as follows. The existence problem is studied in
section 2.  In section 3 the basic properties of the linearized
stability problem are studied and (the decomposition of) the Evans
function is introduced. Section 4 is the main section of the paper, in
it we develop an approach by which the (possible) location and
existence of `slow-fast eigenvalues' near the essential spectrum can
be studied. This section is split into three parts: a subsection on
the regular problem, a subsection in which we study an explicit
example ($G(V) = -\eps^2 \ga$, $H(U^2,V) = H_0 U^2$) in full detail,
and a subsection in which we study the `fate' of the regular front as
$G_1$ approaches $0$ in the general case.  In section 5 we present
simulations which clearly exhibit the impact of the distinction
between a destabilization by the discrete or by the essential
spectrum.  Moreover, we discuss some related issues and topics of
future research.
\begin{Rem} \rm 
\label{rem:moregen}
Large parts of the theory developed in this paper can be generalized
to systems of the type (\ref{most})/(\ref{most_f}) in which the fast
reduced limit system is of the type $U_t = U_{\xi\xi} + B(U^2;V_0)U$
for some function $B$, i.e.  to bi-stable systems of a more general
nature. We focused on he standard case, i.e.  $B = 1 + V_0 - U^2$,
since the analysis is more transparent. If one drops the condition on
the symmetry (\ref{symmU}), the fronts will in general travel with a
certain (nonzero) speed. Although the symmetry is used throughout this
paper, there is no reason to expect that such asymmetric systems
cannot studied along the lines of the methods presented here.
\end{Rem}

\section{The existence problem}\label{sec:exist}
\setcounter{Theorem}{0}
\setcounter{equation}{0}

We analyze the existence of stationary one-dimensional patterns
through geometric singular perturbation theory \cite{fen79,jon95}
using the methods developed in \cite{dkz,dgk3}.  Therefore, we write
the ODE associated to (\ref{decomp}) as a dynamical system in
$\mathbb{R}^4$, \beq
  \label{ODE_f}
  \left\{
  \begin{array}{rcl}
  \dot{u} &=& p
  \\
  \dot{p} &=& -(1+v-u^2)u
  \\
  \dot{v} &=& \eps q
  \\
  \dot{q} &=& \eps \left[-(1+v-u^2)H(u^2,v; \eps) - G(v;\eps) \right],
\end{array}
\right.  
\eeq 
where $\dot { }$ denotes the derivative with respect to the spatial
variable $\xi$ (\ref{fastvar}) (i.e. $\xi$ `plays the role of time').
Note that this system inherits two symmetries of (\ref{decomp})
\beq
\label{ODEsymm}
\xi \to -\xi, p \to -p, q \to -q \; \; {\rm and} \; \; u \to -u, p \to -p. 
\eeq
We consider the `super-slow' case in which $G_1(\eps) = \O(\eps^2)$
separately in sections \ref{sec:ex_suslo} and \ref{sec:gen_suslo}.
Note that in the fast reduced limit, i.e. $\eps\to 0$ in (\ref{ODE_f}), 
the monotonically increasing heteroclinic front solution
is given by $(u_0,p_0,v_0,q_0)$, where
\begin{equation}\label{u0p0}
 (u_0(\xi;v_0),p_0(\xi;v_0)) = \left( \sqrt{1+v_0} \tanh\left(
     \sqrt{\frac{1+v_0}{2} }\xi\right),
 \frac{1+v_0}{\sqrt{2}} \sech^2 \left(
   \sqrt{\frac{1+v_0}{2}}\xi\right) \right)\,,
\end{equation}
and $v_0$ and $q_0$ are constants.

\subsection{The regular case}\label{sec:ex_reg} 
The main result of this section is,
\begin{Theorem}
\label{th:ex_reg}
Let $G_1(\eps)$ (\ref{def:G1H0}) be $\O(1)$ and negative. Then, for
$\eps > 0$ small enough, system (\ref{ODE_f}) has a symmetric pair of
heteroclinic orbits: $\Ga^+_h(\xi;\eps) =
(u_h(\xi;\eps),p_h(\xi;\eps),v_h(\xi;\eps),q_h(\xi;\eps))$ and
$\Ga^-_h(\xi;\eps) =
(-u_h(\xi;\eps),-p_h(\xi;\eps),v_h(\xi;\eps),q_h(\xi;\eps))$, 
with \\
$\lim_{\xi \to \pm \infty}\Ga^+_h(\xi;\eps) = (\pm 1,0,0,0)$ and
$\lim_{\xi \to \pm \infty}\Ga^-_h(\xi;\eps) = (\mp 1,0,0,0)$;
$u_h(\xi;\eps)$ and $q_h(\xi;\eps)$ are odd and monotonic as functions
of $\xi$, $v_h(\xi;\eps)$ and $p_h(\xi;\eps)$ even. Moreover,
$|u_h(\xi;\eps) - u_0(\xi;0)| = \O(\eps)$ (\ref{u0p0})
uniformly on $\mathbb{R}$, $|v_h(\xi;\eps)|, |q_h(\xi;\eps)| =
\O(\eps)$ uniformly on $\mathbb{R}$, and $v_h(0;\eps)$ is the extremal
value of $v_h(\xi;\eps)$, with
\beq 
\label{appr_vheps}
v_h(0;\eps) = \frac{\eps}{2 \sqrt{-G_1(0)}}
\int_{-\infty}^{\infty} \left(1-u^2_0(\xi;0)\right)H(u_0^2(\xi;0),0) d\xi + \O(\eps^2).
\eeq
The orbits $\Ga^{\pm}(\xi;\eps)$ correspond to the (stationary) front
patterns $(\pm U_h(\xi;\eps), V_h(\xi;\eps))$ of (\ref{decomp}) with
$U_h(\xi;\eps)=u_h(\xi;\eps)$ odd as function of $\xi$, $V_h(\xi;\eps)
= v_h(\xi;\eps) = \O(\eps)$ even, $\lim_{\xi \to \pm
  \infty}U_h(\xi;\eps) = \pm 1$, and $\lim_{\xi \to \pm
  \infty}V_h(\xi;\eps) = 0$.
\end{Theorem}
\begin{proof}  
As the system is singularly perturbed, we also consider (\ref{ODE_f})
with the slow scaling $x =\eps \xi$,
(\ref{ODE_f}) is given by,
\beq
\label{ODE_s}
\left\{
\begin{array}{rcl}
\eps u' &=& p 
\\
\eps p' &=& -(1+v-u^2)u 
\\
v' &=&  q
\\
q' &=& \left[-(1+v-u^2)H(u^2,v; \eps) - G(v;\eps) \right],
\end{array}
\right.
\eeq
where $'$ refers to differentiation with respect $x$.
System (\ref{ODE_s}) is referred to as the slow system. 
We begin by finding the locally invariant manifolds of
(\ref{ODE_s}) in the limit $\eps\to 0$.  In this limit, the first two
equations of (\ref{ODE_s}) will reduce to,
\begin{equation}
              p = 0\,,\quad
    -(1+v-u^2)u = 0\,.
\end{equation}
The manifold given by $(u,p,v,q)=(0,0,v,q)$ is not normally
hyperbolic and will not be considered.  However, the manifolds,
denoted $\M^\pm_0$, determined by $(u,p,v,q)=(\pm \sqrt{1+v},0,v,q)$
are normally hyperbolic and thus by \cite{fen79,jon95}, (\ref{ODE_s})
possesses locally invariant manifolds $\M^{\pm}_\eps$, which are
$\O(\eps)$ close to $\M^\pm_0$.  We now determine the leading order
correction to theses manifolds.  Let the manifold $\M^{\pm}_\eps$, be
given by, 
\beq 
\label{Manexan} 
\M^\pm_\eps=\{u=\pm \sqrt{1+v}+\eps
U^\pm(v,q;\eps) ,p=\eps P^\pm(v,q;\eps),v,q\}\,.  
\eeq 
To obtain
successive approximations of $\M_\eps^\pm$, we can expand 
$U^\pm=u_1^\pm+\eps u_2^\pm+\cdots$, and
$P^\pm=p^\pm_1+\eps p_2^\pm+\cdots$.
Using the first two lines of (\ref{ODE_s}) we find,
\begin{equation}\label{u1u2}
p_1^\pm=\frac{q}{2\sqrt{1+v}}, \, p_2^\pm=\frac{\partial
      u_1^\pm}{\partial v}q -\frac{\partial
      u_1^\pm}{\partial q} G(v;\eps), \,
u_1^\pm  =0, \,
u_2^\pm  =\mp\frac{q^2}{4(1+v)^{5/2}}\mp\frac{G(v;\eps)}{(1+v)^{3/2}}.
\end{equation}
Hence, the (slow) flow on the slow manifold is given by,
\begin{equation}
v''=-G(v;\eps)+\O(\eps^2)\,.
\end{equation}
To leading order, this flow is integrable.  The point $(v,q)=(0,0)$, that
corresponds to $(\pm 1, 0, 0, 0)$, is
a critical point on $\M_\eps^\pm$.  Since $G_1<0$, $(0,0)$ is a saddle
on $\M^{\pm}_{\eps}$
with stable direction $(1,\sqrt{-G_1})$ and unstable direction
$(-1,\sqrt{-G_1})$.

A heteroclinic orbit $\Gamma_h^\pm$ from $(\mp 1,0,0,0)$ to
$(\pm1,0,0,0)$ is both an element of $W^u(\M_\eps^\mp)$ and of
$W^s(\M_\eps^\pm)$.  Here we will only consider $\Gamma_h^+$.  The
existence of $\Gamma_h^-$ follows from the symmetry (\ref{ODEsymm}).
The orbit $\Gamma_h^+$ remains exponentially close to $W^u(-1,0,0,0)|_{\M_\eps^-}$
before it 'takes off' and a makes a
'jump' through the fast field.  After that, it 'touches down' on
$\M_\eps^+$ and remains exponentially close to it (and to 
$W^u(1,0,0,0)|_{\M_\eps^+}$ -- see Figure \ref{fig:regto}. The change in $q$
by the passage through the fast field is $\O(\eps)$ (\ref{ODE_f}),
therefore $\Gamma_h^+$ must take off from $\M_\eps^-$ and touch down on 
$\M_\eps^+$ with a $q$-coordinate that is $\O(\eps)$.  Since
$\Gamma_h^+$ is asymptotic to the saddle points
$(0,0)\in\M_\eps^\pm$, it follows that the $v$-coordinate of
$\Gamma_h^+$ must also be $\O(\eps)$.
Note that we have used here implicitly that $G_1 = \O(1)$.

We will determine whether such a trajectory, as $\Gamma_h^+$, is possible
using a Melnikov method.  Both $W^u(\M_\eps^-)$ and
$W^s(\M_\eps^+)$ are $\O(\eps)$ close to the family of heteroclinic
orbits in the fast reduced limit of (\ref{ODE_f}) given in
(\ref{u0p0}). The leading order distance between $W^u(\M_\eps^-)$ and
$W^s(\M^+_\eps)$ can be determined by a Melnikov function for slowly
varying systems \cite{rob}.  
Both $W^u(\M_\eps^-)$ and $W^s(\M^+_\eps)$ intersect
the hyperplane $\{u=0\}$ transversally.  Note that $W^{s,u}(\M_\eps^\pm)\cap\{u=0\}$
is 2-dimensional, thus, since $\{u=0\}$ is 3-dimensional, one expects
a 1-dimensional intersection $W^u(\M_\eps^-)\cap W^s(\M^+_\eps)\cap\{u=0\}$.
The separation between $W^u(\M_\eps^-)$ and
$W^s(\M^+_\eps)$ is, at leading order, measured by the 
integral,
\begin{equation}\label{Mel1}
\Delta = \int_{-\infty}^\infty \left( \begin{array}{c}
p(\xi) \\ u(\xi)+u^3(\xi)-u(\xi)v_0 \end{array} \right) \wedge \left(
\begin{array}{c} 0 \\ -u(\xi)\frac{\partial q}{\partial \delta}(\xi)
\end{array} \right) \,d\xi\,.
\end{equation}
Here the wedge product refers to the scalar cross product and
$\frac{\partial q}{\partial \delta}$ solves the differential equation,
$\frac{d}{d\xi}\left(\frac{\partial q}{\partial \delta}\right)=q_0
\xi\,,\quad \frac{\partial q}{\partial \delta}(0)=0$.  Substituting
(\ref{u0p0}) into (\ref{Mel1}) results in the following expression for the
leading order splitting distance,
\[
\Delta = -\int_{-\infty}^\infty \frac{q
  \xi}{\sqrt{2}}\tanh\left(\frac{\xi}{\sqrt{2}}\right)
\sech^2\left(\frac{\xi}{\sqrt{2}}\right)\,d\xi
       \quad=-q_0\sqrt{2}\,.
\]
Thus, $W^u(\M_\eps^-)\cap W^s(\M_\eps^+) \cap \{u=0\}$ must be
$\O(\eps)$ close to $q=0$.  By the symmetries (\ref{ODEsymm}), we
conclude that $W^u(\M_\eps^-)\cap W^s(\M_\eps^+) \cap \{u=0\}$ must 
be identically $q=0$.  Hence, again by (\ref{ODEsymm}), any solution
that connects $\M_\eps^-$ to $\M_\eps^+$ must have a $u$ component 
that is odd with respect to $\xi$ and a $v$ component that is even
with respect to $\xi$.

We are now ready to determine the take off (touch down) curves
$T_o^-\subset\M^-_\eps$ ($T_d^+\subset \M^+_\eps$) \cite{dkz,dgk3}.  
These curves represent
the points at which the one-dimensional family of orbits 
in $W^u(\M_\eps^-)\cap W^s(\M_\eps^+)$ leave (land
on) $\M_\eps^\pm$. Let the elements of this family be denoted 
$\gamma(\xi;p)$, where the parameter $p > 0$ corresponds to the 
$p$-component of $\gamma(\xi;p)$ as it 
crosses through $\{u=q=0\}$. Note that the
$\gamma$-family forms the Fenichel fibering of 
$W^u(\M_\eps^-) \cap W^s(\M_\eps^+)$ \cite{fen79} and that each
$\gamma(\xi;p)$ is 
asymptotically close to an unperturbed orbit given in
~(\ref{u0p0}). To each $\gamma(\xi;p)$ we associate two orbits,
$\gamma_{\M_\eps^-}(\xi;p)\subset\M^-_\eps$ and
$\gamma_{\M_\eps^+}(\xi;p)\subset\M_\eps^+$ by the fact that
$||\gamma(\xi;p)-\gamma_{\M_\eps^\pm}(\xi;p)||$ is
exponentially small if $\pm \xi>\O(\eps^{-1})$.  We define
$T_o^-$ and $T_d^+$ as the collections of base points of the Fenichel
fibers on $\M_\eps^-$ and on $\M_\eps^-$,
\begin{equation}\label{offdown}
T_o^-=\bigcup_{p>0}  \gamma_{\M_\eps^-}(0;p)\,,\quad 
T_d^+=\bigcup_{p>0}  \gamma_{\M_\eps^+}(0;p)\,.
\end{equation}
We can compute the leading order structure of
$T_o^-$ and $T_d^+$ by considering the effect of the journey through
the fast field on the slow variables $v$ and $q$. Since $v_\xi=\eps q$
and $q=\O(\eps)$ it follows that the change in $v$ through the fast
field is of higher order, i.e. $\O(\eps^2)$.  By construction, $q$
will be an odd function of $\xi$, thus the value of $q$ for a given
$v$ on $T_o^-$ must be $-\frac 1 2 \Delta q(v)$, where $\Delta q(v)$ is the 
change in $q$ due to one full pass through the fast field (during which
$v$ remains (at leading order) constant, $v=v_0$). Similarly, the value of $q$
on $T_d^-$ must be $\frac 1 2 \Delta q(v)$. Since we already know that
both $v$ and $q$ must be $\O(\eps)$ in this regular case, we  
compute $\Delta q(0)$ (by (\ref{ODE_f}), (\ref{u0p0})),  
\[
\Delta q (0) =\int_{-\infty}^\infty \dot{q}|_{v=0}\,d\xi\,,\\
=-\eps \int_{-\infty}^\infty
\left[1-\tanh^2\left(\frac{\xi}{\sqrt{2}}\right)\right]
H(\tanh^2\left(\frac{\xi}{\sqrt{2}}\right),0)\,d\xi+\O(\eps^2)\,.
\]
To establish the existence of the 
heteroclinic orbit $\Ga^+_h(\xi)$, we consider the intersection
$T_{o}^-\cap W^{u}(-1,0,0,0)|_{\M^-_\eps}$ on $\M_{\eps}^-$,
$\O(\eps)$ close to $(-1,0,0,0)$. Thus, $T_{o}^-$ and $W^{u}(-1,0,0,0)|_{\M^-_\eps}$ 
are given by $\{q = -\frac 1 2 \Delta q(0) + \O(\eps^2)\}$ and 
$\{q=\sqrt{-G_1}v+\O(\eps^2)\}$. Figure ~\ref{fig:regto} shows the superposition of $T_o^-$
with $W^u(-1,0,0,0)|_{\M^-_\eps}$ and of $T_d^+$ with
$W^s(1,0,0,0)|_{\M_\eps^+}$. The $v$-coordinate of $T_{o}^-\cap W^{u}(-1,0,0,0)|_{\M^-_\eps}$
is given in (\ref{appr_vheps}).
\begin{figure}[ht]
\label{fig:regto}
\begin{center}
\psfrag{v}{$v$}
\psfrag{q}{$q$}
\psfrag{ws}{$W^{s}(1,0,0,0)|_{\M_\eps^+}$}
\psfrag{wu}{$W^{u}(-1,0,0,0)|_{\M_\eps^-}$}
\psfrag{to}{$T_o^-$}
\psfrag{td}{$T_d^+$}
\psfrag{mp}{$\M_\eps^+$}
\psfrag{mm}{$\M_\eps^-$}
\includegraphics[scale=.5]{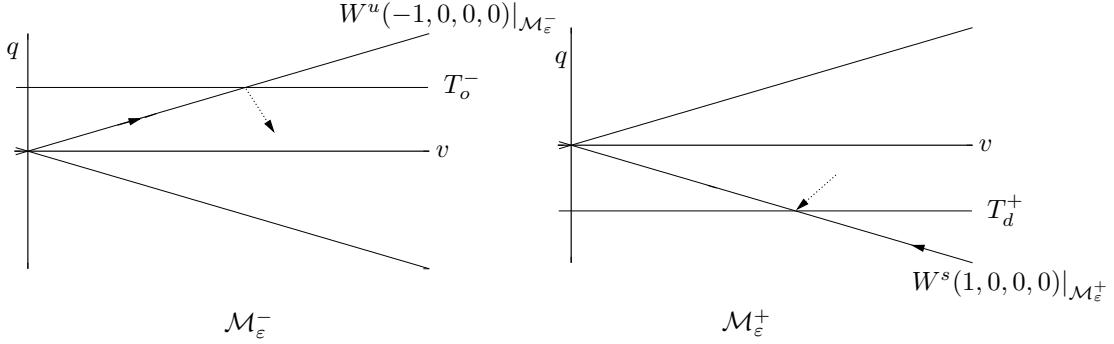}
\caption{Superposition of the 
take off and touch down curves $T_{o, d}^{\pm}$ with 
$W^{s,u}(\pm1,0,0,0)|_{\M_\eps^\pm}$. The intersections
$T_{o}^{-} \cap W^{u}(-1,0,0,0)|_{\M_\eps^-}$
and $T_{d}^{+} \cap W^{s}(1,0,0,0)|_{\M_\eps^+}$
determine the heteroclinic front solution $\Ga_h^+(\xi;\eps)$.  The
dotted arrows inidicates the orbit $\Ga^+_h(\xi)$ 'taking off' and
'touching down'.}
\end{center}
\end{figure}
Thus, we have established the existence of an orbit
$\Gamma^+_h \in W^u(\M_\eps^-)\cap W^s(\M_\eps^+)$ that is asymptotic
to $(-1,0,0,0)\in \M_\eps^-$.  Since $\Gamma_h^+$ passes through
$\{u=0,q=0\}$ during its jump through the fast field, it follows by
the symmetries (\ref{ODEsymm}) that $\Gamma_h^+$ is indeed the orbit
described in the statement of the Theorem.  As was already mentioned,
the existence of $\Gamma_h^-$ also follows immediately from
(\ref{ODEsymm}).  
\hfill\end{proof}\\
\begin{Rem} 
\label{rem:sigdeg} \rm
  We note that if $G_1=\O(\eps^\sigma)$ for some $\sigma\in[0,2)$ then the
  intersection of $T_o^-$ and $W^u(-1,0,0,0)|_{\M^-_\eps}$ will result
  in a value of $v_0$ of $\O(\eps^{2-\sigma}) \ll \O(1)$ (\ref{appr_vheps}), 
  thus $\Ga_h^{+}(\xi)$ will still be a regular perturbation of the
  orbit in the scalar limit. Moreover, this argument also shows that
  singular orbits may exist for $G_1=\O(\eps^2)$.
\end{Rem}
\subsection{The super-slow limit: an example}
\label{sec:ex_suslo}
In this section we consider the `significant degeneration' $G_1(\eps) =
\O(\eps^2)$. For simplicity, we only consider the
case in which the flow on the slow manifolds $\M_\eps^{\pm}$ is linear,
i.e.  $G(v;\eps) = \eps^2 G_1(\eps) \stackrel{\rm def}{=} -\eps^2
\ga$, where $\ga$ does not depend on $\eps$.  Moreover, we first
consider an explicit expression for $H(u^2,v;\eps)$, $H(u^2,v;\eps)=H_0u^2$.  
The case of a general $H(U^2,V)$ will be considered in the next subsection.  
We refer to Remark
\ref{rem:geng1} for a brief discussion of the case of a general
function $G(V)$.  System (\ref{ODE_f}) reduces to 
\beq
\label{ODE_suslo}
\left\{
\begin{array}{rcl}
\dot{u} &=& p
\\
\dot{p} &=& -(1+v-u^2)u
\\
\dot{v} &=& \eps q
\\
\dot{q} &=& \eps \left[-(1+v-u^2)H_0 u^2 + \eps^2 \ga v \right].
\end{array}
\right.
\eeq
This system has various types of (singular) heteroclinic orbits.
\begin{Theorem}
\label{th:ex_suslo}
Assume that $G(V)=-\eps^2\ga V$, $H(U^2,V)=H_0 U^2$ and that $\eps$
is small enough. 
{\bf (i): $H_0>0$.} If $\gamma>\ga_{\rm double}$, where 
$\ga_{\rm double}=\frac{3}{2} H_0^2 + \O(\eps)$, (\ref{ODE_suslo}) has two
pairs of heteroclinic orbits, 
$\Ga^{+,j}_{h}(\xi;\eps) = 
(u_h^{j}(\xi),p_h^{j}(\xi),v_h^{j}(\xi),q_h^{j}(\xi))$, $j = 1,2$,
and their symmetrical counterparts
$\Ga^{-,j}_{h}(\xi;\eps) = 
(-u_h^{j}(\xi),-p_h^{j}(\xi),v_h^{j}(\xi),q_h^{j}(\xi))$,
with $\lim_{\xi \to \pm \infty}\Ga^{+,j}_h(\xi;\eps) = (\pm 1,0,0,0)$.
In the fast field $u_h^{j}(\xi)$, respectively $v_h^{j}(\xi)$, is
asymptotically and uniformly close to $u_0(\xi;v_{j})$
(\ref{u0p0}) resp. $v_{j}$; the constants $v_{j}$ are the zeros of
$\sqrt{\ga} v = \frac23 \sqrt{2} H_0 (v+1)^{3/2}$ so that $0<v_1<2<v_2$ 
(at leading order). In the slow field, $\Ga^{+,j}_{h}(\xi;\eps)$ 
is exponentially close to $W^{u,s}(\pm 1,0,0,0)|_{\M_\eps^\pm} \subset \M_\eps^\pm$.
The orbits $\Ga^{\pm,1}_{h}(\xi;\eps)$ and $\Ga^{\pm, 2}_{h}(\xi;\eps)$
merge in a saddle-node bifurcation of heteroclinic orbits as  
$\gamma \downarrow \ga_{\rm double}$. There are no heteroclinic orbits
for $\gamma<\ga_{\rm double}$.
\\
{\bf (ii): $H_0 < 0$.} The relation $\sqrt{\ga} v = \frac23 \sqrt{2} H_0 (v+1)^{3/2}$
has a unique zero for all $\ga > 0$ and there is one pair of heteroclinic
orbits $\Ga^\pm_h(\xi; \eps)$ for all $\ga > 0$. These orbits have
the same structure as described in {\bf (i)}.
\\
The orbits $\Ga^{\pm(,j)}_h(\xi; \eps)$ correspond to the front solutions
$(U^{\pm(,j)}_h(\xi; \eps), V^{\pm(,j)}_h(\xi; \eps))$ of (\ref{decomp})
with $U^{\pm(,j)}_h(\xi; \eps) = \pm u_h^{j}(\xi;\eps)$ odd, and
$V^{\pm(,j)}_h(\xi; \eps) = v_h^{j}(\xi;\eps)$ even as function of $\xi$.
\end{Theorem}
\begin{proof}  The essence of the analysis of the super-slow
system is similar to that of the regular case.  The important
difference being that, although the change in $q$ by a `jump' 
through the fast field is still $\O(\eps)$, the $v$-coordinate of
the heteroclinic orbit may now be $\O(1)$, due to the super-slow
character of the flow on $\M_{\eps}^{\pm}$.  It is this
difference that will cause the bifurcation and the formation of the
second orbit in case {\bf (i)}.  The flow on the slow manifold is now
$\O(\eps^2)$, i.e. super-slow, and is at leading order governed by,
\begin{equation}
v''=\eps^2 \gamma v\,.
\end{equation}
Since the right hand side of this equation is $\O(\eps^2)$, one might
expect that one needs to incorporate the higher order corrections to
the approximation of $\M^{\pm}_{\eps}$ (\ref{u1u2}) to determine
the leading order flow on $\M^{\pm}_{\eps}$. However, the
$\O(\eps^2)$ correction contains a term with a $q^2$ factor and a term with
$G(v)$ (\ref{u1u2}). Since we consider $q =\O(\eps)$ on $\M^{\pm}_{\eps}$
and since $G(v) = \O(\eps^2)$, the 
resulting correction will not be of leading order.

Again the equilibria on $\M^{\pm}_{\eps}$ are saddles, with stable and unstable
directions, $(\pm1,\eps\sqrt{\gamma})$.  As in Theorem \ref{th:ex_reg} 
we only consider the orbit that jumps from $\M_\eps^-$ to $\M_\eps^+$ 
(the others follows from the symmetry (\ref{ODEsymm}). 
We repeat the Melnikov calculations and again conclude that,
$W^u(\M_\eps^-)\cap W^s(\M_\eps^+) \cap \{u=0\}$ must be
identically $q=0$.  Hence, again by (\ref{ODEsymm}), any solution
that connects $\M_\eps^-$ to $\M_\eps^+$ must have a $u$ component 
that is odd with respect to $\xi$ and a $v$ component that is even
with respect to $\xi$.

We define the take off, $T_o^-$, and touch down, $T_d^+$, curves as in
(\ref{offdown}).  We find the leading order behavior of $T_d^+$ and
$T_o^-$, by calculating the change in $q$ as we traverse the fast
field.  As in the regular case, $v$ remains a constant up to $\O(\eps^2)$ 
and the value of $q$ on the take off (touch down) curve must
be $-\frac 1 2 \Delta q(v_0)$ ($\frac 1 2 \Delta q(v_0)$), where
$v_0$ is the (leading order) constant value of the $v$-coordinate
of the orbit that is heteroclinic to $\M_\eps^+$ in the fast field. 
The calculation of
the change in $q$ is similar to that of the regular case except
that $v_0$ now effects the leading order term (\ref{u0p0}),
\[
\begin{array}{lcl}
\Delta q(v_0) & = & -\eps H_0(1+v_0)^2 \int_{-\infty}^\infty
   \left[1-\tanh^2\left(\sqrt{\frac{v_0+1}{2}}\xi\right)\right]
   \tanh^2\left(\sqrt{\frac{v_0+1}{2}}\xi\right)\,d\xi+\O(\eps^2)\\
        & = &-\eps\frac{2\sqrt{2}}{3}H_0(v_0+1)^{3/2}+\O(\eps^2)\,.
\end{array}
\]
\begin{figure}[ht]\label{fig:ssto}
\begin{center}
\psfrag{k}{$\M_\eps^-$}
\psfrag{rel}{$q$}
\psfrag{com}{$v$}
\subfigure[]
{\includegraphics[scale=.35]{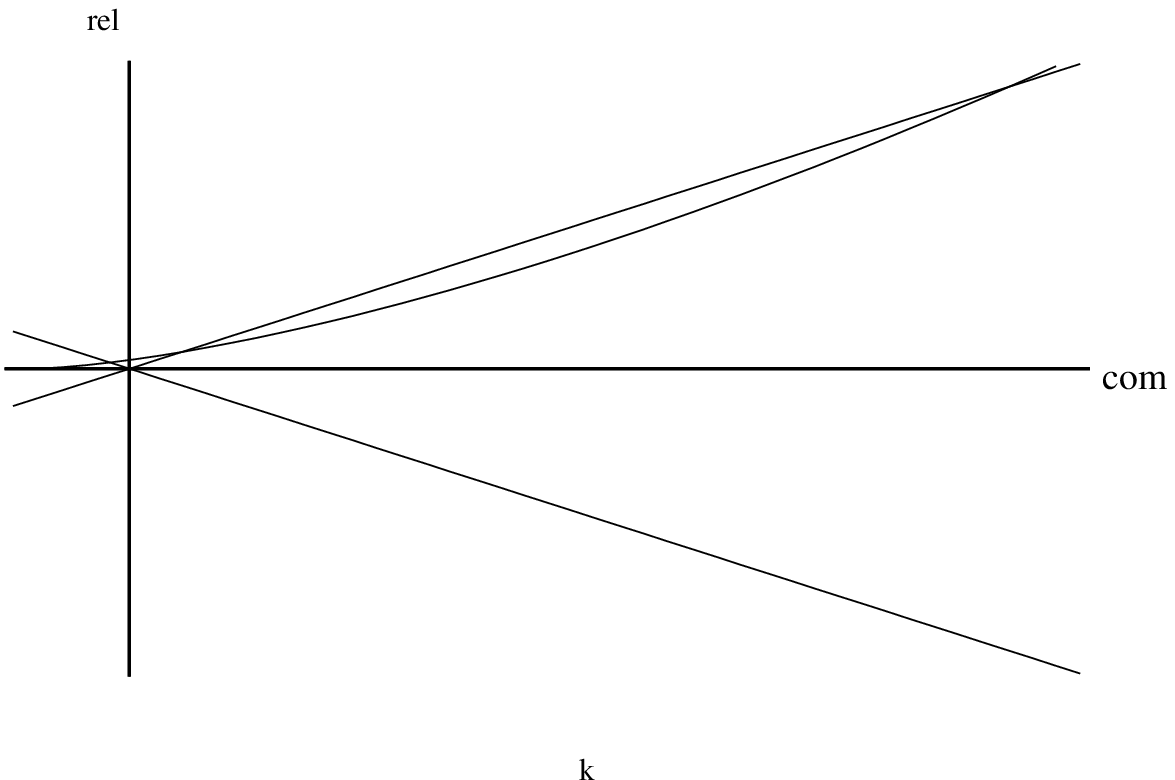}}
\hspace{2cm}
\subfigure[]
{\includegraphics[scale=.35]{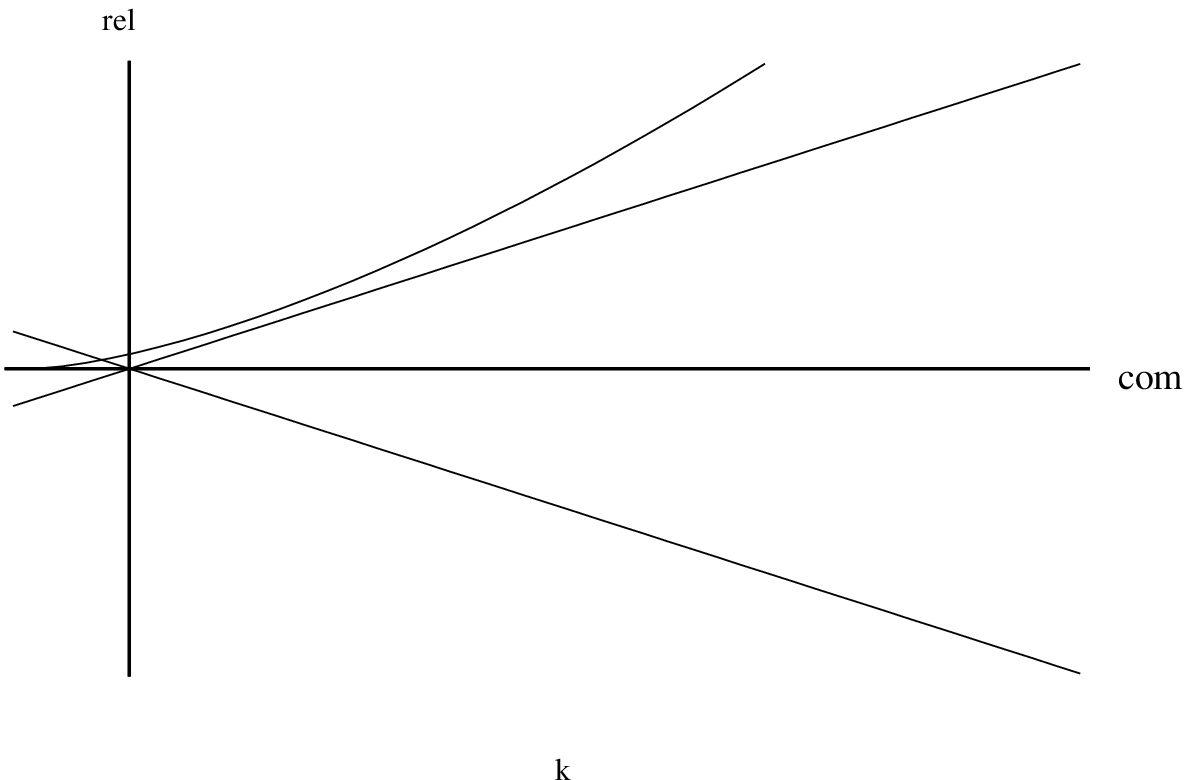}}
\caption{Superposition of $T_o^-$ with $W^{u}(-1,0,0,0)|_{\M_\eps^-}$ in
  the super-slow case with
  $H_0 > 0$, for $\ga > \ga_{\rm double}$ (a) and  $\ga > \ga_{\rm double}$ (b).}
\end{center}
\end{figure}
The heteroclinic orbits are again determined by 
$T_o^- \cap W^{u}(-1,0,0,0)|_{\M^-_\eps}$, where
$T_o^- = \{q = -\frac12 \Delta q(v_0) + \O(\eps^2)\}$ and 
$W^{u}(-1,0,0,0)|_{\M^-_\eps}= \{q=\eps\sqrt{\ga}v +\O(\eps^2)\}$
\begin{equation}\label{v0rel}
\frac{2\sqrt{2}}{3}H_0(v_0+1)^{3/2}=\sqrt{\gamma}v_0\,,
\end{equation}
see Figure ~\ref{fig:ssto}.
Thus, in the super-slow case, a heteroclinic orbit may leave $\M_\eps^-$ 
with a $v$-coordinate of $\O(1)$.
Now if $H_0 > 0$ and 
$\gamma>\gamma_{\rm double} = \frac 3 2 H_0^2 + \O(\eps)$, (\ref{v0rel})
has two possible solutions, $v_0 = v_{j}$, $j = 1,2$, with $0 < v_1 < 2 < v_2$ 
(at leading order). These intersection correspond to the 
heteroclinic orbits $\Ga_h^{+,j}(\xi)$. 
For $\gamma < \gamma_{\rm double}$, there are no solutions to
(\ref{v0rel}) and thus no heteroclinic connections exist:
the orbits $\Ga_h^{+,1}(\xi)$ and $\Ga_h^{+,2}(\xi)$ have coalesced at 
$\gamma = \gamma_{\rm double}$.  In the case 
that $H_0<0$, (\ref{v0rel}) has a unique solution for 
all values of $\ga > 0$, there is only one pair
of heteroclinic orbits.
\hfill\end{proof}\\

\begin{Rem} \rm
\label{rem:geng1}
If $G(V)$ is not linear in the singular limit (i.e. $G_1 = \O(\eps^2)$),
then the analysis becomes more involved, but there are no essentially
new phenomena. In this case, the magnitude (w.r.t. $\eps$) of the second 
derivative of $G(v)$ at $v = 0$  will start to play a role
comparable to $G_1$. Moreover, the flow on $\M_{\eps}^{\pm}$ is nonlinear,
so that $W^{u, s}(\pm 1,0,0,0)|_{\M^-_\eps}$ is no longer a straight line
(at leading order), therefore, many `new' intersections of
$T_o^- \cap W^{u}(-1,0,0,0)|_{\M^-_\eps}$, and thus `new' heteroclinic
orbits, may appear.  
\end{Rem}  
\subsection{The super-slow limit: the general case}\label{sec:gen_suslo}
We now consider the general super-slow problem, i.e. (\ref{ODE_f})
with $G=-\eps^2\ga v$.  The treatment of the general super-slow case
and (\ref{ODE_suslo}) is in essence identical to that of the previous
section. However, the statement of the main results cannot be formulated
as explicit as in Theorem \ref{th:ex_suslo}, as long as there is no
explicit expression given for $H(U^2,V)$. Nevertheless, the character
of the existence result is similar to that of Theorem \ref{th:ex_suslo},
there can be various kinds of heteroclinic orbits that might coalesce
in saddle-node bifurcations.

As in the proofs of Theorems \ref{th:ex_reg} and \ref{th:ex_suslo}, the
existence of the heteroclinic orbits is established by the intersection
of $T_o^-$ and $W^u(-1,0,0,0)|_{\M_{\eps}^-}$, i.e. by the solution
$v_0$ of
\beq
\label{defhets_gen}
\sqrt{\ga} v_0 = \frac12 \int_{-\infty}^{\infty}
\left[1+v_0-u_0^2(\xi; v_0) \right] H(u_0^2(\xi; v_0),v_0) d\xi,
\eeq
at leading order. Note that the right hand side equals $-\frac12 \Delta q(v_0)$,
i.e. half the accumulated change in $q$ during a circuit through the fast field, and
that we have used (\ref{u0p0}). 
\begin{Theorem}
\label{th:gen_suslo}
Assume that $G(V)=-\eps^2\ga V$, and that $\eps$ is small enough.
System (\ref{ODE_f}) has $n\ge 0$ pairs of heteroclinic orbits,
$\Ga^{\pm,j}_{h}(\xi;\eps) = (\pm u_h^{\pm,j}(\xi), \pm
p_h^{\pm,j}(\xi),v_h^{\pm,j}(\xi),q_h^{\pm,j}(\xi))$, where
$j=1,\ldots\,,n$, with $\lim_{\xi \to \pm \infty}\Ga^{+,j}_h(\xi;\eps)
= (\pm 1,0,0,0)$.  The number $n = n(\ga)$ is given by the number of
solutions $v_j$ of (\ref{defhets_gen}).  In the fast field
$u_h^{j}(\xi)$, respectively $v_h^{j}(\xi)$, is asymptotically and
uniformly close to $u_0(\xi;v_{j})$ (\ref{u0p0}) resp. $v_{j}$, where
the constant $v_{j}$ is the $j$-th zero of (\ref{defhets_gen}).  In
the slow field, $\Ga^{+,j}_{h}(\xi;\eps)$ is exponentially close to
$W^{u,s}(\pm 1,0,0,0)|_{\M_\eps^\pm} \subset \M_\eps^\pm$.
\\
Two orbits $\Ga^{\pm,j}_{h}(\xi;\eps)$ and $\Ga^{\pm,
  j+1}_{h}(\xi;\eps)$ coalesce in a saddle-node bifurcation of
heteroclinic orbits at a certain value $\gamma = \ga^j_{\rm double}$,
if the zeroes $v_j \leq v_{j+1}$ of (\ref{defhets_gen}) merge, i.e. if
the intersection $T_o^- \cap W^u(-1,0,0,0)|_{\M_{\eps}^-}$ is
non-transversal.

The orbits $\Ga^{\pm,j}_h(\xi; \eps)$ correspond to
the front solutions $(U^{\pm,j}_h(\xi; \eps), V^{\pm,j}_h(\xi; \eps))$
of (\ref{decomp}) with \\ $U^{\pm,j}_h(\xi; \eps) = \pm
u_h^{j}(\xi;\eps)$ odd, and $V^{\pm,j}_h(\xi; \eps) =
v_h^{j}(\xi;\eps)$ even as function of $\xi$.
\end{Theorem}

\begin{figure}[ht]\label{fig:gensuslo}
\begin{center}
\psfrag{v}{$v$}
\psfrag{q}{$q$}
\subfigure{\includegraphics[scale=.5]{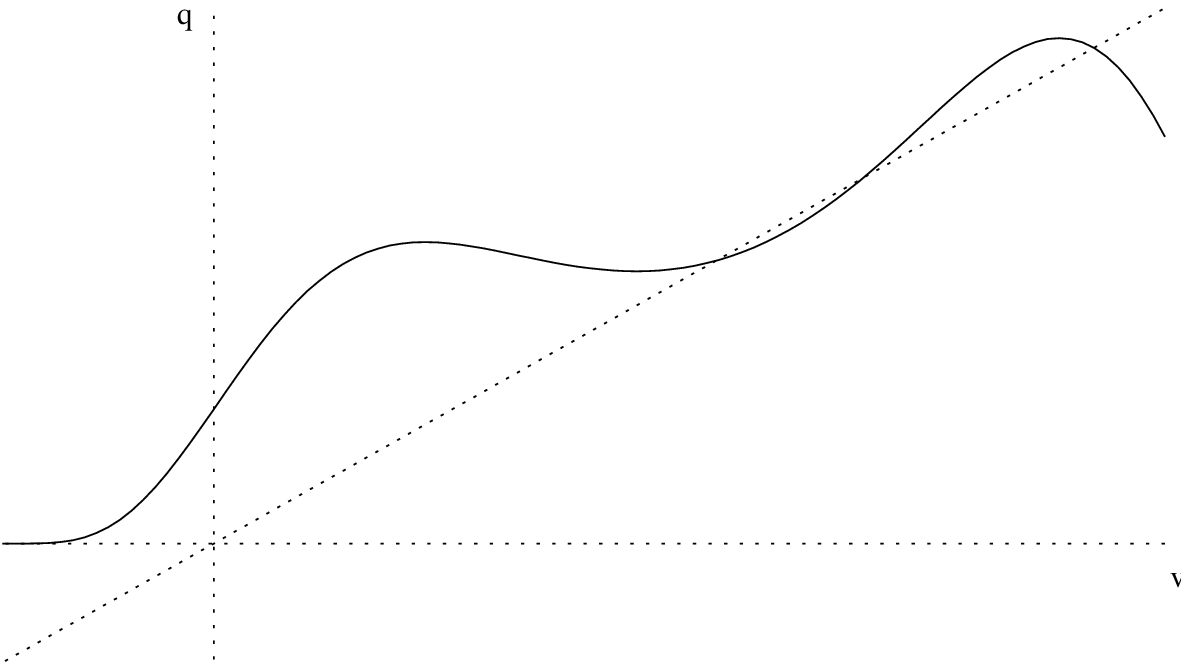}}\hspace{.7cm}
\subfigure{\includegraphics[scale=.5]{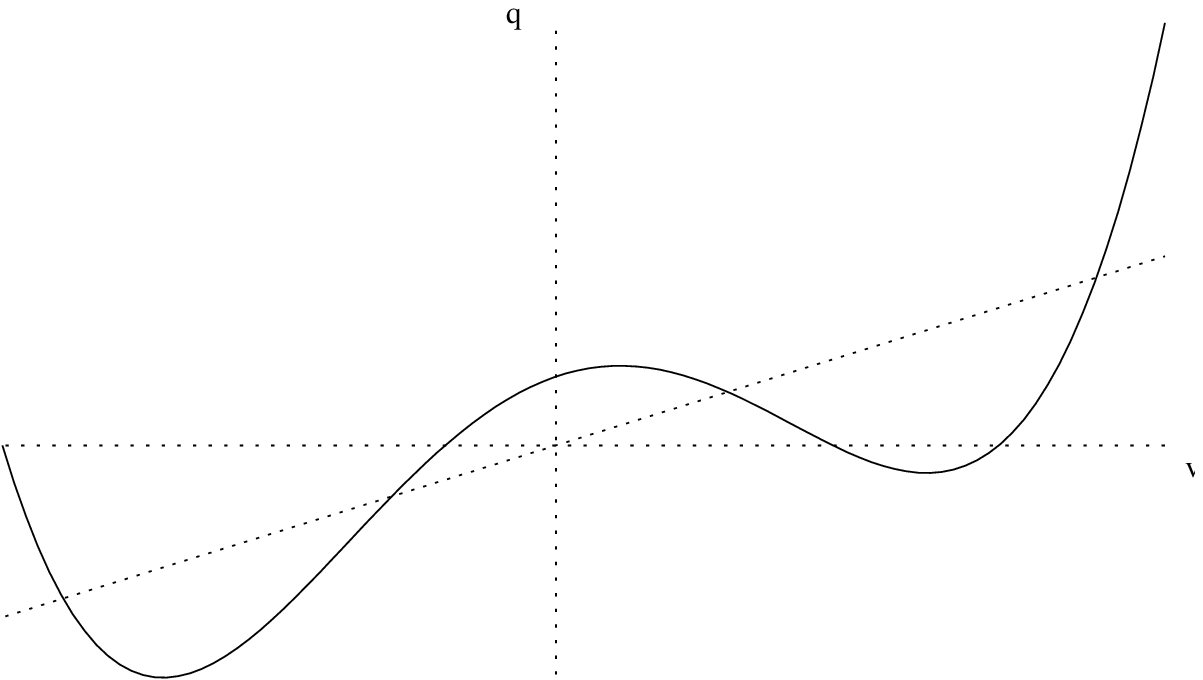}}\vspace{.25cm}
\caption{Two examples of the possible character of the intersection 
  $T_o^- \cap W^{u}(-1,0,0,0)|_{\M_\eps^-}$ for a given $H(U^2,V)$;
  (a) there are 3 different singular heteroclinic orbits (b) 4
  heteroclinic orbits.}
\end{center}
\end{figure}

The proof of this result is in essence identical to that of Theorem \ref{th:ex_suslo}.
In Figure \ref{fig:gensuslo} two examples of the possible richness of
the intersection $T_o^- \cap W^u(-1,0,0,0)|_{\M_{\eps}^-}$ are given.

\section{The stability of fronts}
\label{sec:stab_1}
\setcounter{Theorem}{0} \setcounter{equation}{0}

With a slight abuse of notation we (re-)introduce $u(\xi)$ and $v(\xi)$ by
\[
U(\xi,t) = U_h(\xi;\eps) + u(\xi) e^{\la t}, \;
V(\xi,t) = V_h(\xi;\eps) + v(\xi) e^{\la t},
\]
substitute this into (\ref{decomp}), and linearize
\beq
\label{linstab}
\begin{array}{ccl}
u_{\xi\xi} & + & (1+V_h - 3U_h^2 - \la) u = - U_h v  \\
v_{\xi\xi} & = & \eps^2 \left\{
2 \left[H(U_h^2,V_h) - (1+V_h-U_h^2)
\frac{\pa H}{\pa U^2}(U_h^2,V_h)\right] U_h u \right. \\
& & \left.\, \; \; \; \; - \left[H(U_h^2,V_h) - (1+V_h-U_h^2)
\frac{\pa H}{\pa V}(U_h^2,V_h) + \frac{\pa G}{\pa V}(V_h)
- \tau \la \right] v \right\}.
\end{array}
\eeq
Note that the front pattern $(U_h(\xi), V_h(\xi))$ corresponds to any
of the regular or singular heteroclinic orbits $\Ga^{\pm, j}_h(\xi)$
of Theorems \ref{th:ex_reg}, \ref{th:ex_suslo}, \ref{th:gen_suslo}.
In the stability analysis of forthcoming sections we will only
consider the front patterns of $+$-type, i.e. those fronts for which
$\lim_{\xi \to \pm \infty} U_h(\xi;\eps) = \pm 1$. Thus, we do not 
explicitly consider their symmetric counterparts. 
Due to the symmetry (\ref{symmU}) this
is of course also not necessary. 
The coupled system of second order equations (\ref{linstab})
is equivalent to a linear
system in $\mathbb{C}^4$,
\beq
\label{lindyn}
\phi_{\xi} = A(\xi; \la, \eps) \phi \; \; {\rm with} \; \;
\phi(\xi) = (u(\xi),p(\xi),v(\xi),q(\xi)),
\eeq
where $A(\xi; \la, \eps)$ is a $4 \times 4$ matrix with ${\rm Tr}(A(\xi; \la, \eps) \equiv 0$,
and $u_{\xi}=p$, $v_{\xi} = \eps q$.
It follows that
\beq
\label{Apminfty}
\lim_{\xi \to \pm \infty} A(\xi; \la,\eps) \stackrel{\rm def}{=} A^{\pm}_{\infty}(\la, \eps)
=
\left(
\begin{array}{cccc}
0 &  1 &  0 &  0 \\
2 + \la &  0 &  \mp 1 & 0 \\
0 &  0 & 0 &  \eps \\
\pm 2 \eps H_0 &  0 & -\eps(H_0 + G_1 -\tau \la) &  0
\end{array}
\right) \eeq (\ref{def:G1H0}). The matrices $A^{\pm}_{\infty}$ have
the same set of eigenvalues $\La_i(\la, \eps)$, $i =1,2,3,4$, \beq
\label{Lambdas}
\La_{1,4}^2(\la,\eps) = \la + 2 + \O(\eps^2), \;
\La_{2,3}^2(\la,\eps) = \eps^2 \frac{\tau \la^2 - \la(G_1+H_0-2\tau) - 2G_1}{\la + 2} + \O(\eps^4).
\eeq
Note that both expansions break down as $\la$ approaches $-2$ (see Remark \ref{rem:la=-2}).
We define a branch cut such that for $z \in \mathbb{C}$
${\rm arg}(\sqrt{z}) \in (-\frac12 \pi, \frac12 \pi]$, so that
the $\La_i$'s can be ordered
\beq
\label{Laorder}
{\rm Re}(\La_4(\la, \eps)) < {\rm Re}(\La_3(\la, \eps)) < 0
< {\rm Re}(\La_2(\la, \eps)) < {\rm Re}(\La_3(\la, \eps))
\eeq
This ordering of course breaks down if $\la \in \si_{\rm ess}$, 
the essential spectrum associated to (\ref{linstab})/(\ref{lindyn}), since
$\si_{\rm ess}$ coincides with values of $\la$ for which either
${\rm Re}(\La_{1,4}(\la, \eps)) = 0$ or ${\rm Re}(\La_{2,3}(\la, \eps)) = 0$
\cite{henry}; see also section \ref{sec:stab_ess}).
The eigenvectors $E^{\pm}_i(\eps,\la)$ of the matrices $A^{\pm}_{\infty}(\eps,\la)$
associated to $\La_i(\la, \eps)$ are given by
\beq
\label{Eipm}
E^{\pm}_{1,4}(\eps,\la)
\left(
\begin{array}{c}
1 \\ \La_{1,4}(\la, \eps) \\ \O(\eps^2) \\ \pm \frac{2 H_0}{\La_{1,4}(\la, \eps)} \eps + \O(\eps^3)
\end{array}
\right),
\; \;
E^{\pm}_{2,3}(\eps,\la)
\left(
\begin{array}{c}
\pm \frac{1}{\la + 2} + \O(\eps^2) \\ \O(\eps^2) \\ 1 \\ \frac{1}{\eps} \La_{2,3}(\la, \eps)
\end{array}
\right)
\eeq
(for $\la + 2 = \O(1)$ -- Remark \ref{rem:la=-2}).

\begin{Rem} \rm
\label{rem:la=-2}
The expansions (\ref{Lambdas}) and (\ref{Eipm}) are
only valid for $\la + 2 = \O(1)$ with respect to $\eps$.
It is straightforward to check that $\La_{1,4}^2(\la,\eps) = \O(\eps) = \La_{2,3}^2(\la,\eps)$
if $\la + 2 = \O(\eps)$ and that in general, when $\la + 2 = \O(\eps^{\si})$ for some
$\si \in [0,1]$, $\La_{1,4}^2(\la,\eps) = \O(\eps^{\si})$ and
$\La_{2,3}^2(\la,\eps) = \O(\eps^{2-\si})$. Thus, $\La_{1,4}$ cannot be assumed to be
large/fast compared to $\La_{2,3}$ if $\la + 2 = \O(\eps)$. Since $\la \approx -2$ is
way into the stable region, we do not consider this degeneration further and assume
throughout this paper that $|\La_{2,3}| \ll |\La_{1,4}|$.
\end{Rem}

\subsection{The essential spectrum}
\label{sec:stab_ess}

\begin{figure}[ht]
\label{fig:essspec}
\begin{center}
\psfrag{k}{}
\psfrag{rel}{}
\psfrag{com}{}
\subfigure[Re$(\lambda)$ vs $k$ with $H_0<0$.]
{\includegraphics[scale=.35]{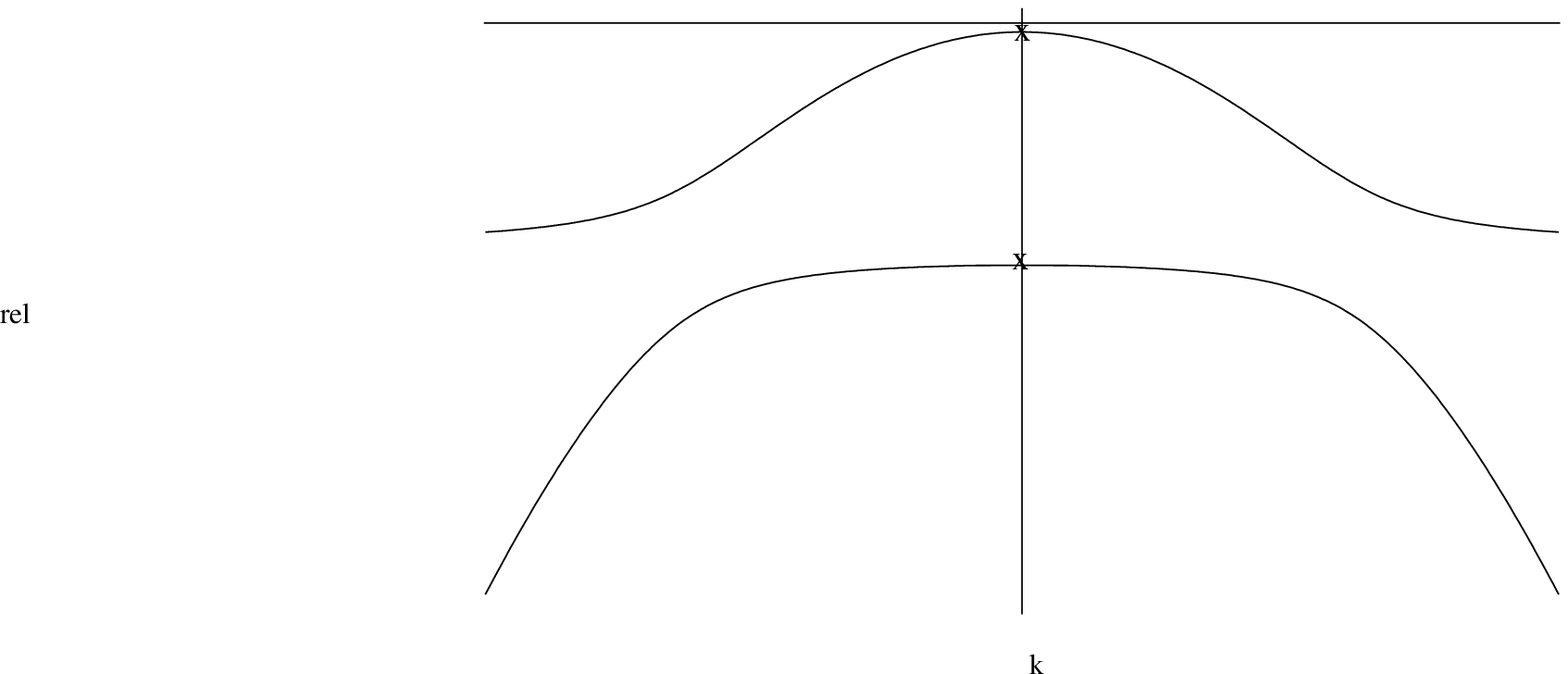}}
\hspace{2cm}
\subfigure[Re$(\lam)$ vs Im$(\lam)$ with $H_0<0$.]
{\includegraphics[scale=.35]{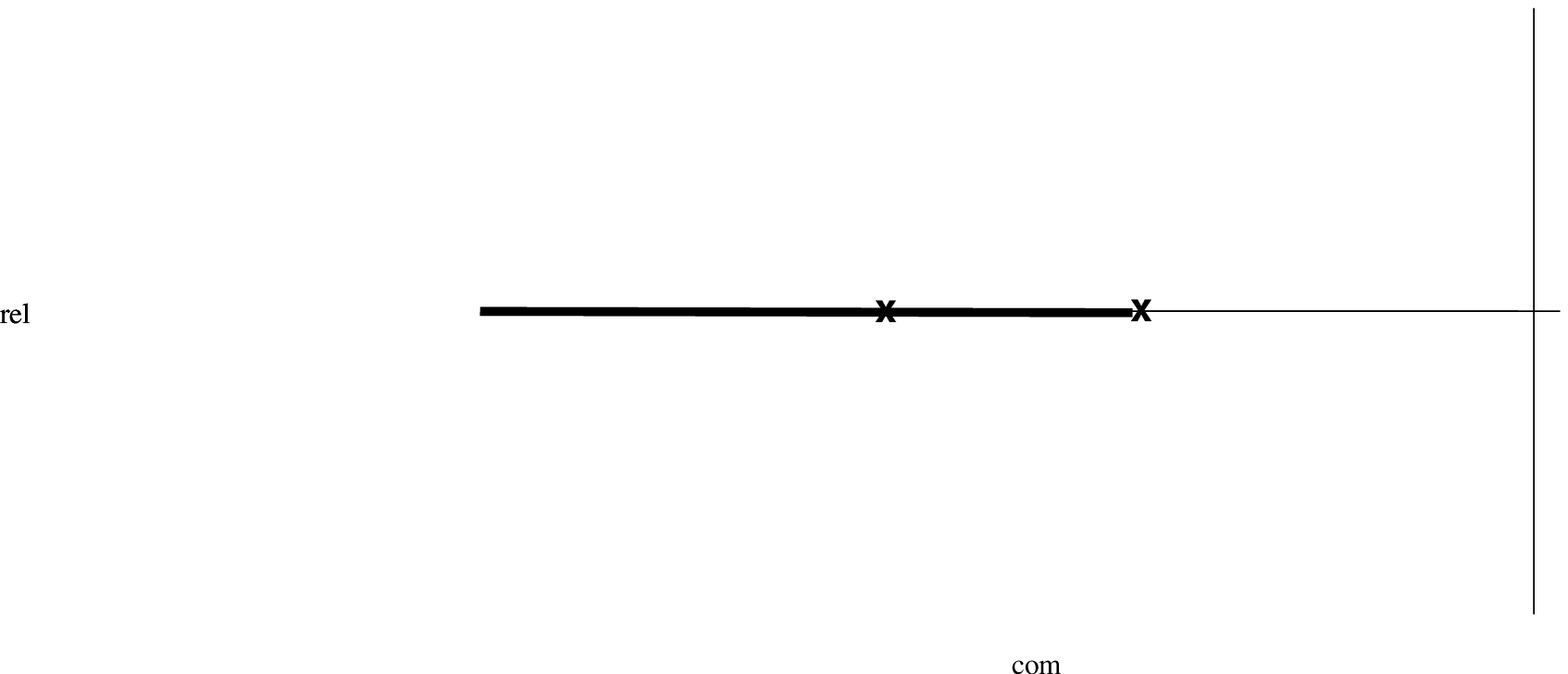}}
\subfigure[$H_0=0$.]
{\includegraphics[scale=.35]{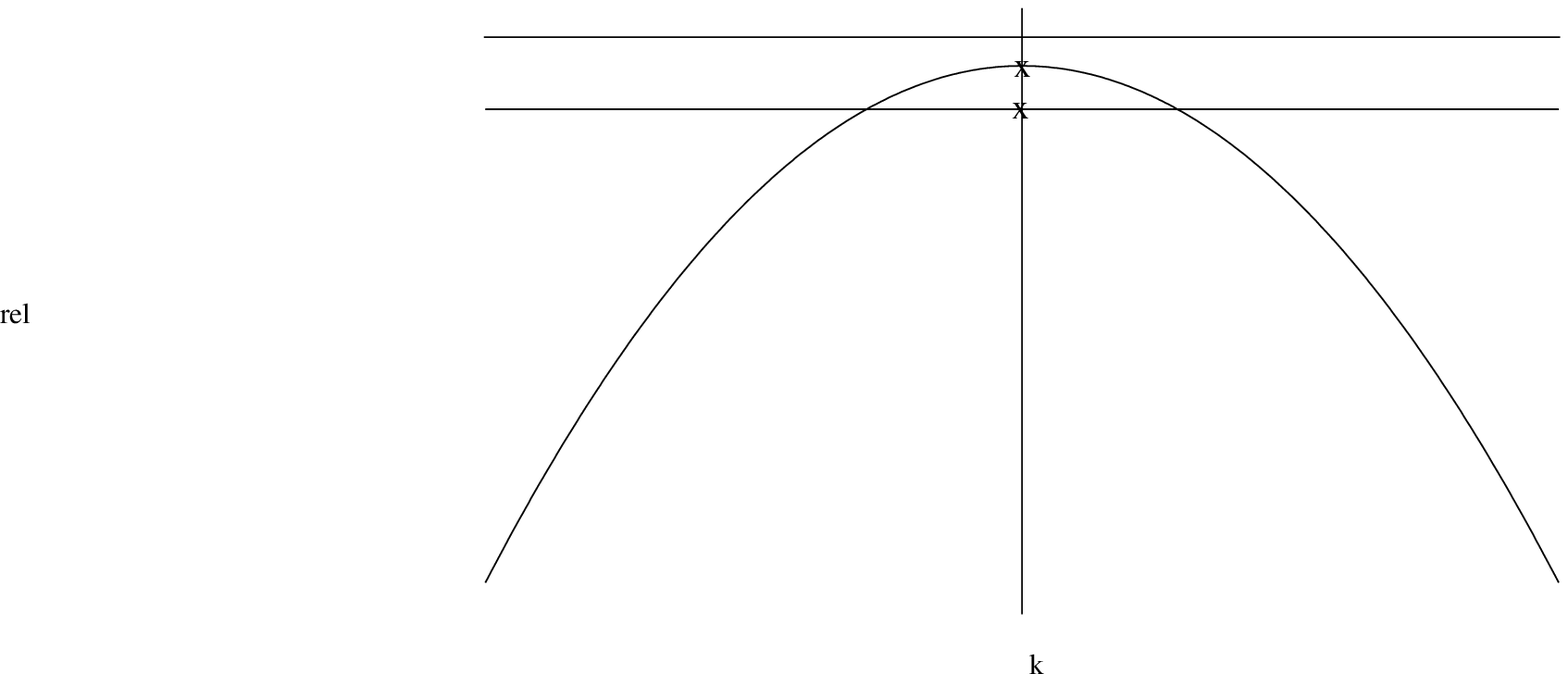}
\hspace{2cm}}
\subfigure[$H_0=0$.]
{\includegraphics[scale=.35]{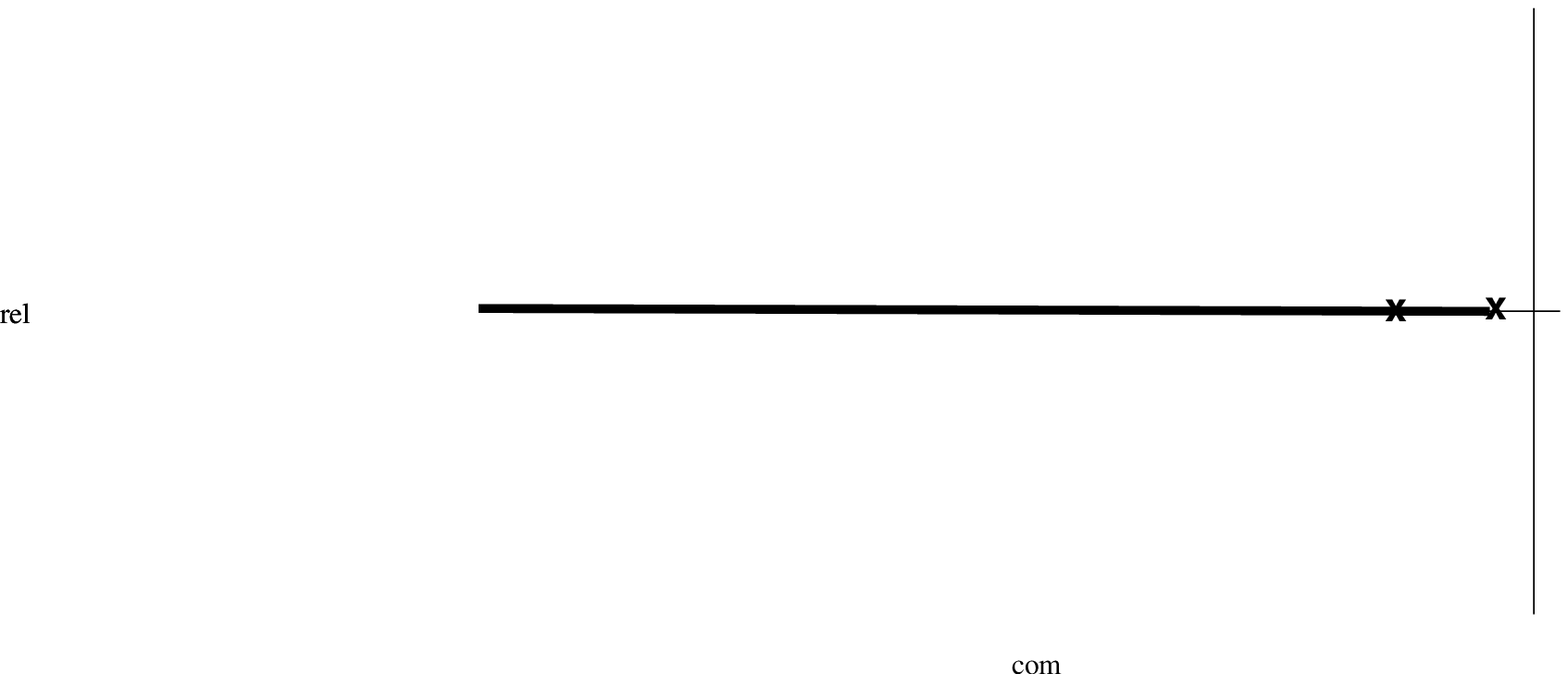}}
\subfigure[$H_0 \in (0,(\sqrt{2\tau} - \sqrt{-G_1})^2)$.]
{\includegraphics[scale=.35]{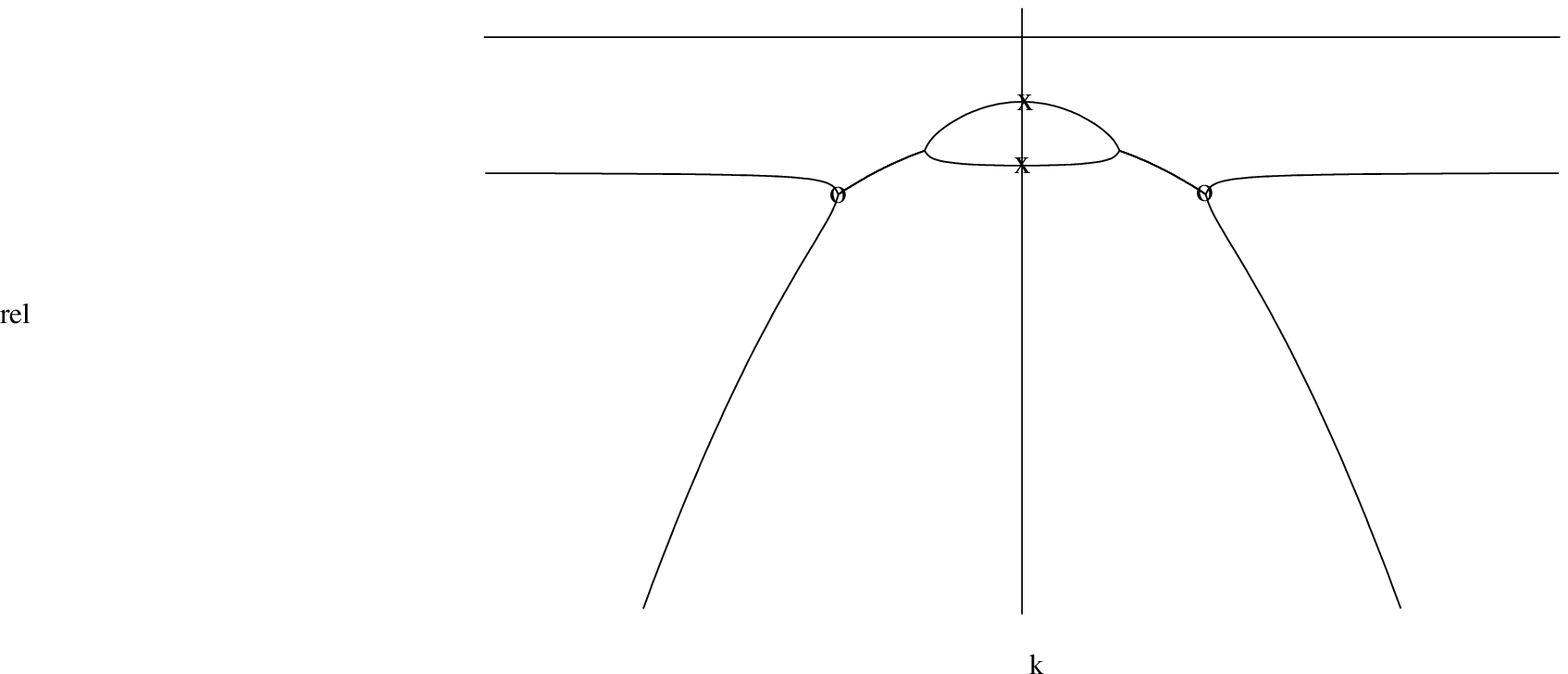}
\hspace{2cm}}
\subfigure[$H_0 \in (0,(\sqrt{2\tau} - \sqrt{-G_1})^2)$.]
{\includegraphics[scale=.35]{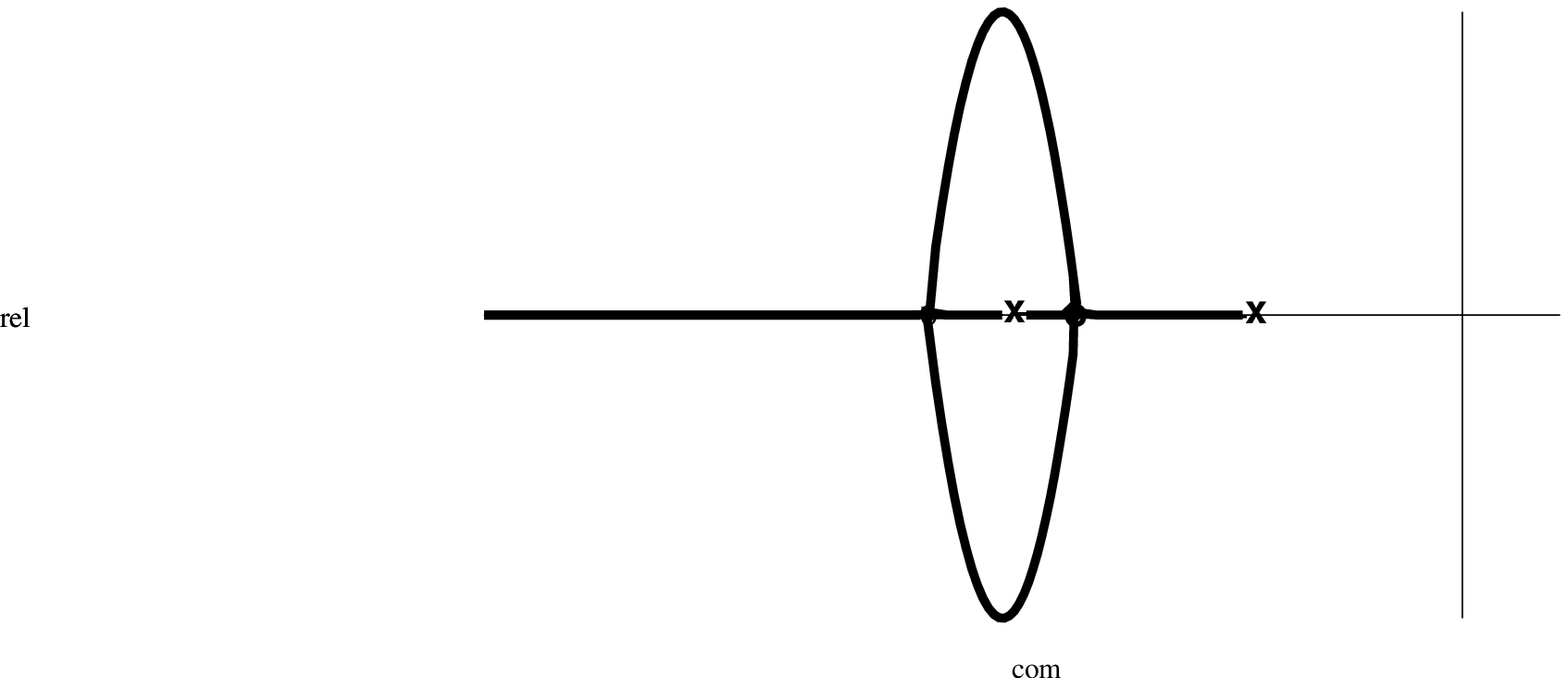}}
\subfigure[$H_0=(\sqrt{2\tau} - \sqrt{-G_1})^2$.]
{\includegraphics[scale=.35]{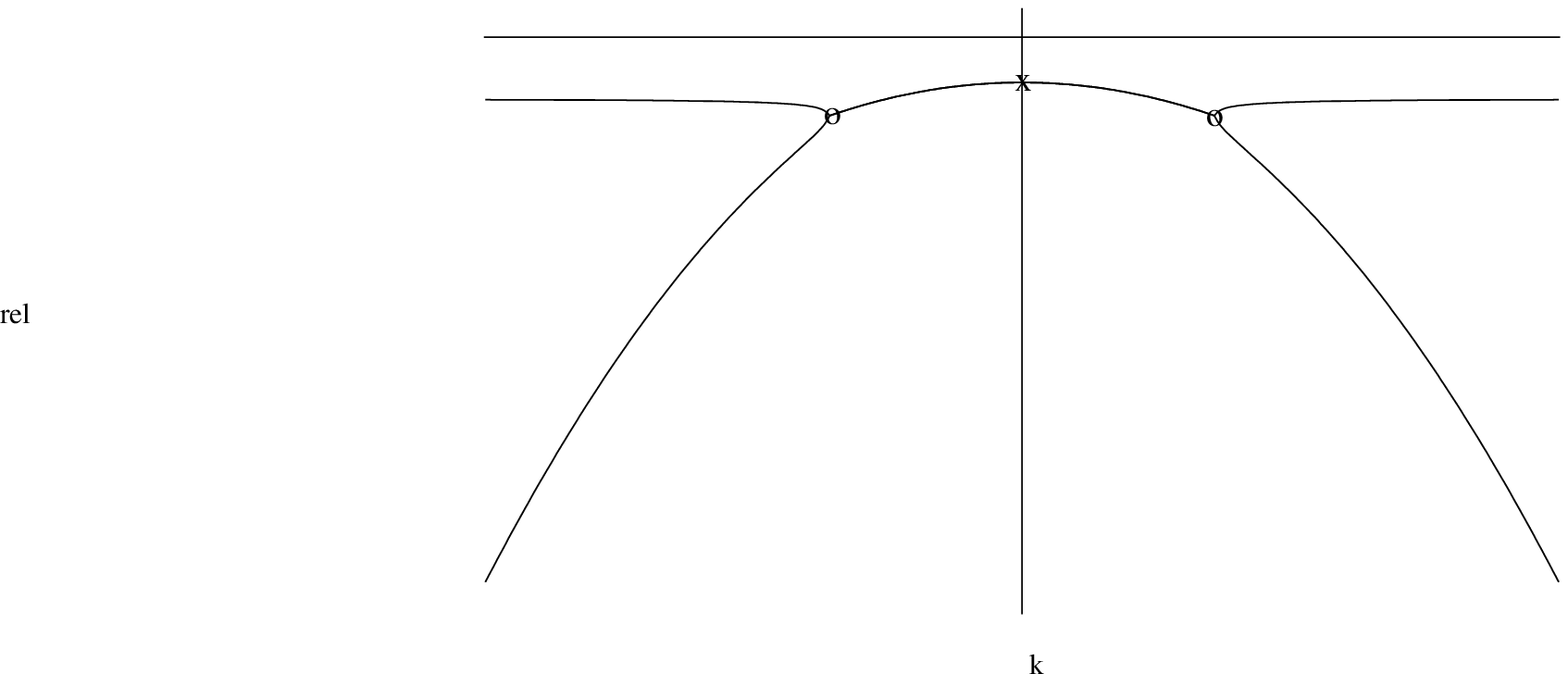}
\hspace{2cm}}
\subfigure[$H_0=(\sqrt{2\tau} - \sqrt{-G_1})^2$.]
{\includegraphics[scale=.35]{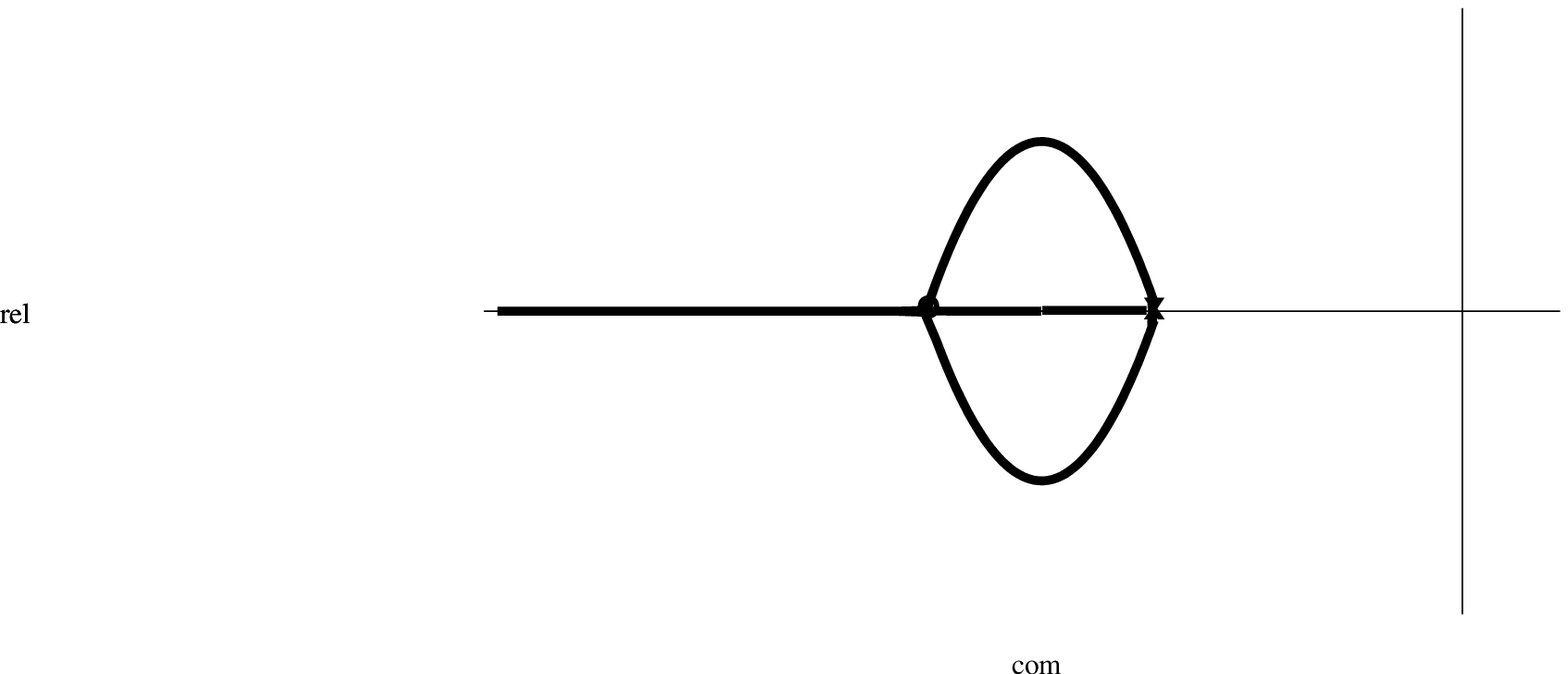}}
\subfigure[$H_0 \in ((\sqrt{2\tau} - \sqrt{-G_1})^2, 2\tau - G_1)$.]
{\includegraphics[scale=.35]{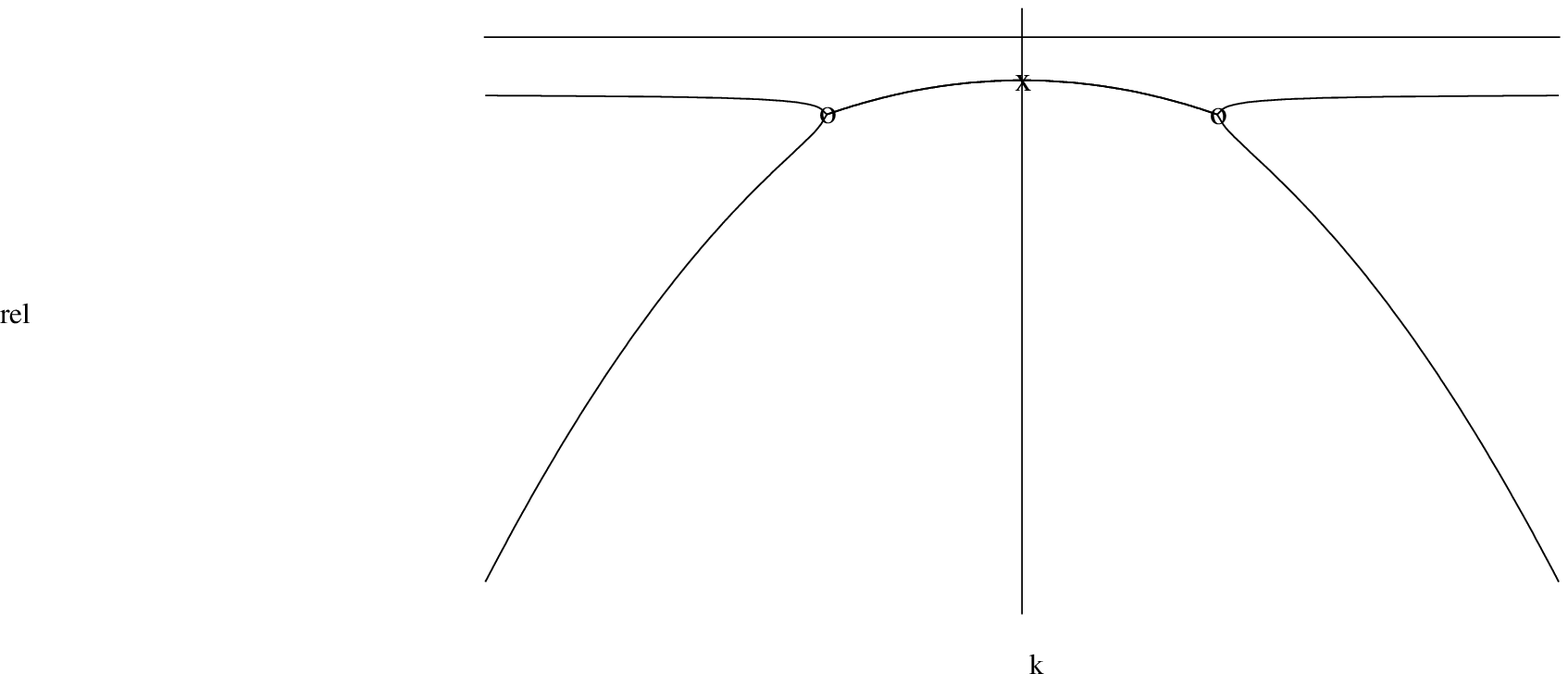}
\hspace{2cm}}
\subfigure[$H_0 \in ((\sqrt{2\tau} - \sqrt{-G_1})^2, 2\tau - G_1)$.]
{\includegraphics[scale=.35]{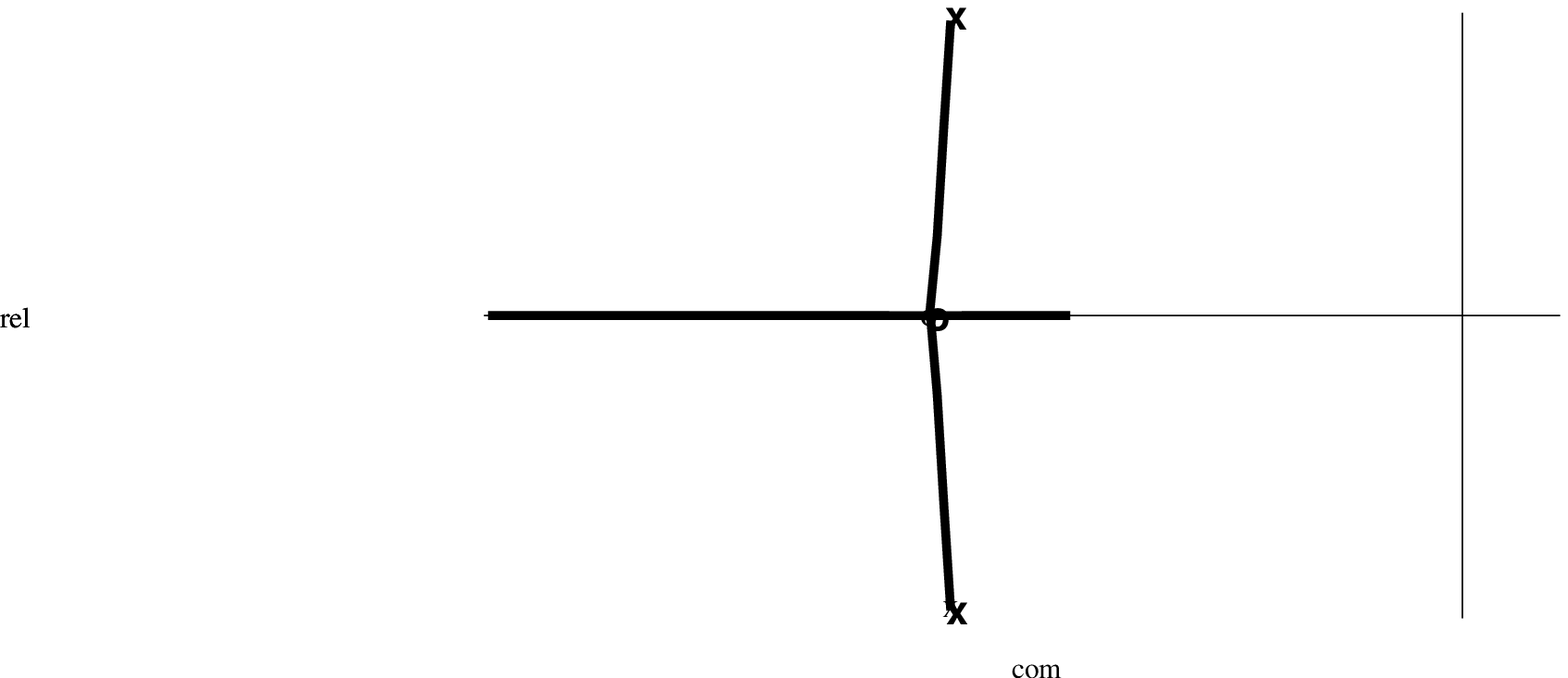}}
\caption{The five possible different structures of the stable essential spectrum.  On the left we plot Re$(\la)$ vs $k$ and on the right Re$(\la)$ vs Im$(\la)$.}
\end{center}
\end{figure}

The essential spectrum associated to the stability of the front patterns
$(U,V) = (U_h(\xi), V_h(\xi))$ is fully determined by the spectrum of
the linear stability problem for
the (trivial) background states (at $\pm \infty$) $(U,V) \equiv (\pm 1,0)$
\cite{henry}. Therefore, we introduce $k \in \mathbb{R}$ and
$\al, \be, \la \in \mathbb{C}$ by
\[
U(x,t) = \pm 1 + \al e^{ik\xi + \la t}, \; \; V(x,t) = \be e^{ik\xi + \la t},
\]
and substitute this expression into (\ref{decomp}) (using (\ref{fastvar})).
This yields the matrix equation
\[
\begin{pmatrix}
-k^2 - 2& \pm 1\\
\mp 2 \eps^2 H_0 & - k^2 + \eps^2(H_0 + G_1)\\
\end{pmatrix}
\begin{pmatrix}
\al \\ \be
\end{pmatrix}
= \la
\begin{pmatrix}
\al \\ \eps^2 \tau \be
\end{pmatrix},
\]
where
$G_1$ and $H_0$ have been introduced in (\ref{def:G1H0}). Thus,
$\la = \la(k^2)$ is a solution of the characteristic equation
\beq
\label{chareq}
Q(\la,k) =
(\la+k^2+2)(\eps^2\tau\la+k^2-\eps^2 (H_0+G_1)) + 2\eps^2H_0
= 0.
\eeq
Note that this equation holds for both background states $(\pm 1,0)$,
due to the symmetry (\ref{symmU}).
We may conclude
\begin{Lem}
\label{lem:stabsig_ess}
The essential spectrum $\si_{\rm ess}$ associated to (\ref{linstab}) is
given by the solutions $\la = \la(k^2)$ of (\ref{chareq}) with $k \in \mathbb{R}$;
$\si_{\rm ess}$ is stable, i.e. $\si_{\rm ess} \in \{{\rm Re}(\la) < 0\}$,
if $G_1 < 0$ and $H_0+G_1-2\tau < 0$.
\end{Lem}
\begin{proof}
The two conditions in this lemma are obtained directly from
\beq
\label{la1la2}
\begin{array}{cccccl}
\la_1 + \la_2 & = & \frac{1}{\eps^2\tau} 
\left[\eps^2(H_0 + G_1 - 2\tau)-k^2(1+\eps^2 \tau)\right] & < & 0 &
\forall \; k
\\
\la_1 \la_2 & = & \frac{1}{\eps^2\tau} 
\left[k^4 + k^2(2-\eps^2(H_0+G_1)) - 2\eps^2 G_1\right] & > & 0 &
\forall \; k.
\end{array}
\eeq Both relations attain their extremal value at $k=0$. 
\hfill\end{proof}\\
However, we need to have more information on the essential spectrum
than just this stability result. In section \ref{sec:stab_slfa} we
will see that the appearance of edge bifurcations is closely related
to the structure of $\si_{\rm ess}$. We focus on the stable case $G_1 < 0$ and
$H_0 + G_1 - 2 \tau < 0$. It is straightforward to check
that (\ref{chareq}) has two solution $\la_{1,2}(k) \in \mathbb{R}$
for all $k \in \mathbb{R}$ if $H_0 < 0$. As $H_0$ passes through
zero two $k$-intervals, $(-k^+, -k^-)$ and $(k^-, k^+)$ ($0 < k^- < k^+$)
appear in which $\la_{1,2}(k)$ are complex valued. These intervals
merge (i.e. $k^- \downarrow 0$) as $H_0$ approaches 
$(\sqrt{2 \tau} - \sqrt{-G_1})^2$. For 
$(\sqrt{2 \tau} - \sqrt{-G_1})^2 < H_0 < 2 \tau - G_1$ (which
is a non-empty region),  $\la_{1,2}(k) \in \mathbb{C}$
if $-k^+ < k < k^+$. See Figure \ref{fig:essspec}.

\subsection{The Evans function}
\label{sec:Evans}
The use of the Evans function in the analysis of linear systems
associated to the stability of traveling waves is by now
well-established. Here, we give a brief exposition of the
characteristics of the Evans function in reaction-diffusion systems.
We refer to \cite{agj,pw,gj,dgk2,dgk3} for the full analytic details
of the statements in this section.
\\ \\
We define the complement of the essential spectrum by \beq
\label{def:Ce}
\C_e = \mathbb{C} \backslash \si_{\rm ess}.
\eeq
For $\la \in \C_e$ the ordering (\ref{Laorder}) holds, so that
\begin{Lem}
\label{lem:decomp}
For all $\la \in \C_e$
there exist two two-dimensional families of solutions
$\Phi_-(\xi; \la, \eps)$ and $\Phi_+(\xi; \la, \eps)$
to (\ref{lindyn}) such that
$\lim_{\xi \to \pm \infty} \phi_{\pm} (\xi; \la, \eps) = (0,0,0,0)^t$
for all $\phi_{\pm} (\xi; \la, \eps) \in \Phi_{\pm}(\xi; \la, \eps)$;
$\Phi_{\pm}(\xi; \la, \eps)$ depend analytically on $\la$.
\end{Lem}

An eigenfunction of (\ref{lindyn}) must be in the intersection of
$\Phi_{-}(\xi; \la, \eps)$ and $\Phi_{+}(\xi; \la, \eps)$. Since Tr$(A) \equiv 0$,
we therefore define
the Evans function $\D(\la, \eps)$ by
\beq
\label{def:D}
\D(\la, \eps) = {\rm det}[\phi_1(\xi; \la, \eps), \phi_2(\xi; \la, \eps),
\phi_3(\xi; \la, \eps), \phi_4(\xi; \la, \eps)],
\eeq
where $\{\phi_{1}, \phi_{2}\}$ (respectively
$\{\phi_{3}, \phi_{4}\}$) span the space
$\Phi_-$ (resp. $\Phi_+$). The
Evans function is analytic in $\la \in \C_e$, its zeroes
correspond to eigenvalues of (\ref{lindyn}) counting multiplicities \cite{agj,pw}.
Of course, this definition does not determine $\D(\la)$ uniquely.
However, this can be achieved by choosing $\phi_{1}(\xi)$ and $\phi_{2}(\xi)$
as follows
\begin{Lem}
\label{lem:defphi12t12}
For all $\la \in \C_e$
there is a unique solution
$\phi_1(\xi; \la, \eps) \in \Phi_- (\xi; \la, \eps)$ of (\ref{lindyn})
such that
\[
\lim_{\xi \to -\infty} \phi_1(\xi; \la, \eps) e^{-\La_1(\la, \eps) \xi} =
E^-_1(\la, \eps).
\]
(\ref{Lambdas}), (\ref{Eipm}). There exists an analytic transmission function
$t_1(\la, \eps)$ such that
\[
\lim_{\xi \to \infty} \phi_1(\xi; \la, \eps) e^{-\La_1(\la, \eps) \xi} =
t_1(\la, \eps) E^+_1(\la, \eps).
\]
For $\la \in \C_e$ such that $t_1(\la, \eps) \neq 0$
there is a unique solution
$\phi_2(\xi; \la, \eps) \in \Phi_- (\xi; \la, \eps)$ of (\ref{lindyn}),
that is independent of $\phi_1(\xi; \la, \eps)$, that satisfies
\[
\lim_{\xi \to -\infty} \phi_2(\xi; \la, \eps) e^{-\La_2(\la, \eps) \xi} =
E^-_2(\la, \eps)
\; \; {\rm and} \; \; 
\lim_{\xi \to \infty} \phi_2(\xi; \la, \eps) e^{-\La_1(\la, \eps) \xi} =
(0,0,0,0)^t.
\]
There exists a second meromorphic transmission function
$t_2(\la, \eps)$, that is determined by
\[
\lim_{\xi \to \infty} \phi_2(\xi; \la, \eps) e^{-\La_2(\la, \eps) \xi} =
t_2(\la, \eps) E^+_2(\la, \eps).
\]
\end{Lem}
The solutions $\phi_{3,4}(\xi; \la, \eps) \in \Phi_+ (\xi; \la, \eps)$ of (\ref{lindyn})
can be defined likewise. Since 
$\sum_{i=1}^{4} \La_i(\la, \eps) \equiv 0$ (\ref{Lambdas}),
\[
\D(\la, \eps)= {\rm det}
[\phi_1(\xi) e^{-\La_1 \xi}, 
\phi_2(\xi) e^{-\La_2 \xi}, \phi_3(\xi) e^{-\La_3 \xi}, 
\phi_4(\xi) e^{-\La_4 \xi}], 
\]
so that $\D(\la, \eps)$ 
can be decomposed into a product
of $t_1(\la, \eps)$ and $t_2(\la, \eps)$ by taking the limit $\xi \to +\infty$.
\begin{Lem}
\label{lem:decD}
Let $\la \in \C_e$, then
\beq
\label{decD}
\D(\la, \eps) = t_1(\la, \eps) t_2(\la, \eps) \, {\rm det}
\left[ E^+_1(\la, \eps), E^+_2(\la, \eps), E^+_3(\la, \eps), E^+_4(\la, \eps) \right].
\eeq
\end{Lem}
We conclude that the eigenvalues of (\ref{lindyn}) correspond to zeroes of the
transmission functions $t_1(\la, \eps)$ and $t_2(\la, \eps)$. However, we will
see that a zero of $t_1(\la, \eps)$ does not necessarily correspond to a
zero of $\D(\la, \eps)$, since $t_2(\la, \eps)$ can have poles (see also \cite{dgk2,dgk3}).

\subsection{The fast eigenvalues}
\label{sec:stab_fast}
Next section will be devoted to the analysis of (the zeroes of)
$t_2(\la, \eps)$, here we consider the zeroes of the fast transmission
function $t_1(\la, \eps)$.  In order to do so, we first consider the
stability problem associated to the front solution $U_f(\xi; V_0)$,
with $U_f(\xi; V_0) \to \pm (1 + V_0)$ as $\xi \to \pm \infty$, of the scalar
fast reduced limit problem (\ref{scalarlim}), 
\beq
\label{linfast}
w_{\xi\xi} + (1 + V_0 - 3u_0^2(\xi; V_0) - \la) w = 0,
\eeq
since $U_f(\xi; V_0) = u_0(\xi; V_0)$ (\ref{u0p0}).
This system can be written as a linear system in $\mathbb{C}^2$,
\beq
\label{lindyn_f}
\psi_{\xi} = B(\xi; \la) \psi \; \; {\rm with} \; \; \psi(\xi) =
(u(\xi),p(\xi)), \eeq where $B(\xi; \la)$ is a $2 \times 2$ matrix of
which the coefficients are by construction $\O(\eps)$ close (uniformly
in $\xi$) to those of the $2 \times 2$ block in the upper left corner
of the $4 \times 4$ matrix $A^{\pm}(\xi; \la, \eps)$ defined in
(\ref{lindyn}), if we set $V_0 = V_h(0)$. The Evans function
associated to this problem can be written as $\D_f(\la) = {\rm
  det}[\psi_1(\xi, \la), \psi_4(\xi, \la)]$, in which $\psi_1(\xi)$
and $\psi_4(\xi)$ are solutions of (\ref{lindyn}) determined by
$\lim_{\xi \to -\infty} \psi_1(\xi) e^{-\sqrt{\la + 2}\xi} = (1,
\sqrt{\la + 2})^t$ and $\lim_{\xi \to \infty} \psi_4(\xi)
e^{\sqrt{\la+2}\xi} = (1, -\sqrt{\la+2})^t$ (where $\pm \sqrt{\la +
  2}$ and $(1, \pm \sqrt{\la + 2})^t$ are the eigenvalues and
eigenvectors of the matrix $B_{\infty}(\la) = \lim_{\xi \to \pm
  \infty} B(\xi; \la)$ (compare to (\ref{Lambdas}), (\ref{Eipm}))). As
for the full system, we can define an analytic fast reduced
transmission function $t_f(\la)$ by $\lim_{\xi \to \infty} \psi_1(\xi)
e^{-\sqrt{\la + 2}\xi} = t_f(\la)(1, \sqrt{\la + 2})^t$, so that
\[
\D_f(\la) =  \lim_{\xi \to \infty} {\rm det}[\psi_1(\xi), \psi_4(\xi)]
= {\rm det} [t_f(\la) (1, \sqrt{\la+2})^t, (1, -\sqrt{\la+2})^t]
= -2 t_f(\la) \sqrt{\la+2}.
\]
The transmission function $t_1(\la)$ is, by construction, asymptotically
close to its fast reduced limit $t_f(\la)$.
\begin{Lem}
\label{lem:zeroest1}
Let $\la_i^f \in \C_e$ such that $t_f(\la_i^f) = 0$. There is a uniquely determined
$\la_i(\eps)$ with $\lim_{\eps \to 0} \la_i(\eps) = \la_i^f$ such that 
$t_1(\la_i(\eps), \eps) = 0$;
$t_1(\la, \eps) \neq 0$ for $\la \neq \la_i(\eps)$.
\end{Lem}

The proof of this Lemma is completely analogous to the proofs of similar statements in
\cite{agj,gj,dgk2,dgk3}.
Hence, we find (the leading order behavior of) the zeroes of $t_1(\la, \eps)$ by
computing the spectrum of (\ref{linfast}). By (\ref{u0p0}) and by introducing
$\eta = \sqrt{\frac12 (1+V_0)}$ we can write (\ref{linfast}) as
\[
w_{\eta \eta} + \left( \frac{6}{\cosh^2 \eta} - P^2 \right) w = 0 \; \; {\rm with} \; \;
P^2 = \frac{2 \la}{1+V_0} + 4,
\]
which is a well-studied problem of Schr\"odinger/Sturm-Liouville type
(see for instance \cite{titch, dgk3}). It has discrete spectrum at $P=1$ and $P=4$
and essential spectrum for $P \in i \mathbb{R}$. We conclude that the eigenvalues of
(\ref{linfast}), and thus the leading order approximations of the zeroes of $t_1(\la)$,
are given by
\beq
\label{fastevs}
\la^f_1 = 0, \; \; \la^f_2 = -\frac{3}{2} (1+V_0) < 0.
\eeq
The essential spectrum of (\ref{linfast}) is given by
\beq
\label{siessf}
\si_{\rm ess}^f = \{ \la \leq -2 (1+V_0) \}.
\eeq
We conclude this subsection by stating two simple, but useful results:
\begin{Lem}
\label{lem:evenodd}
Let $(u(\xi;\eps),v(\xi;\eps))$ be a pair of eigenfunction solutions of (\ref{linstab})
associated to a simple eigenvalue $\la(\eps)$, then either $u(\xi)$
is even as function of $\xi$ and $v(\xi)$ odd, or $u(\xi)$ is odd and
$v(\xi)$ even.
\end{Lem}
\begin{proof}
We write (\ref{linstab}) in the following way,
\begin{gather}\label{recall}
v_{\xi\xi} =\eps^2[F_o(\xi) u + F_e(\xi) v]\,.
\end{gather}
By construction, $U_h$ is an odd function of $\xi$ and $V_h$ is an
even function of $\xi$.  It thus follows that the above functions,
$F_o$ and $F_e$ must be odd and even functions of $\xi$ respectively.
Let $(u,v)$ be an eigenfunctions associated to the eigenvalue $\lam$.
We decompose $(u,v)$ into odd and even components, $u=u_o+u_e$,
$v=v_o+v_e$, where $u_o$, $v_o$ are odd and $u_e$, $v_e$ are even.  By
the parity of the functions $U_h$, $V_h$, $F_o$ and $F_e$ it is clear
that $(u_o,v_e)$ and $(u_e,v_o)$ form two independent solutions of the
eigenvalue problem associated to the eigenvalue $\la$.  Since we have
assumed that $\lam$ is simple, we have a contradiction.
\hfill\end{proof}\\
\begin{Lem}
\label{lem:lais0}
Assume that the eigenfunction solution $v(\xi)$ of (\ref{linstab})
with eigenvalue $\la(\eps)$ is odd, then $\la(\eps) \equiv 0$, so that
$(u(\xi),v(\xi)) = (U_{h,\xi}(\xi;\eps), V_{h,\xi}(\xi;\eps))$.  
\end{Lem}

We will see in section \ref{sec:stab_slfa} that there can be several
eigenvalues for which $u(\xi)$ is odd and $v(\xi)$ even.
\begin{proof}
It is clear that there is an eigenvalue $\la = 0$ associated to the
derivative of the front $(u(\xi),v(\xi)) = (U_{h,\xi}(\xi;\eps),
V_{h,\xi}(\xi;\eps))$. We assume there is another eigenfunction with
$v$ odd.  Since $v_{\xi \xi}$ is $\O(\eps^2)$ and $v$ is odd, it
follows that $|v|\ll 1$ on the fast spatial scale.  Hence, the
equation for the $u$-component is to leading order homogeneous and
given by (\ref{linfast}) (with $u$ replaced by $w$). 
Lemma \ref{lem:evenodd} implies that $u$
is even.  Since the only even eigenfunction of (\ref{linfast}) is
$U_{h,\xi}$ with eigenvalue $0$, it follows that the leading order
behavior of $u$ is given by $U_{h,\xi}$ and that $\lam$ is
asymptotically close to $0$.  We thus write,
\begin{equation}\label{uexpan2}
u=U_{h,\xi}+\delta(\eps)u_1\,,\quad
v=\delta(\eps)v_1\,,\quad
\lam=\delta(\eps)\hlam(\epsilon)\,,
\end{equation}
where $\delta(\eps)\to 0$ as $\eps\to 0$ and $\hlam(0)\ne 0$ (i.e.
$\delta(\eps)$ represents the leading order magnitude of $\hlam$).  We
substitute (\ref{uexpan2}) into (\ref{linstab}), to get the following
equation for $u_1$,
\[
u_{1,\xi \xi}+(1-3U_h^2)u_1=\hlam U_{h,\xi}-U_h v_1\,.
\]
This equation has the solvability condition,
$\int_{-\infty}^{\infty} ( \hlam U_{h,\xi}-U_h v_1 )U_{h,\xi}\,dy=0$.
Now, $v_1$ and $U_h$ are odd while $U_{h,\xi}$ is even, thus we have
that $\hlam=0$, contradicting our assumption.  So the only possible
eigenfunctions with $v$ odd must correspond to a $0$ eigenvalue and
hence, $(u,v)=(U_{h,\xi},V_{h,\xi})$.
\hfill\end{proof}\\

\section{Slow-fast eigenvalues and edge bifurcations}
\label{sec:stab_slfa}
\setcounter{Theorem}{0}
\setcounter{equation}{0}

The `slow-fast eigenvalues' are the eigenvalues that exist due to
the interaction of the fast $U$-equation and the slow $V$-equation in (\ref{decomp}),
thus, these eigenvalues do not have a counterpart in the fast reduced scalar limit 
problem (\ref{scalarlim}). The slow-fast eigenvalues correspond to the zeroes of the
$t_2(\la;\eps)$, since this transmission function is based on a balance between slow
and fast effects. See also Remark \ref{rem:eigenfcts}.

In order to study the combined effect of slow and fast terms, we need to define
the region in which the fast $\xi$-jump takes place more accurately 
\beq
\label{def:If}
I_f = \{\xi \in (-\frac{1}{\sqrt{\eps}}, \frac{1}{\sqrt{\eps}}) \} \; \;
{\rm or} \; \;  \{x \in (-\sqrt{\eps}, \sqrt{\eps}) \}
\eeq
(\ref{fastvar}). Note that the exact choice of the boundaries of $I_f$
is not relevant, any choice will be suitable as long as it is in the transition
zone between $x$ and $\xi$ (i.e. on the boundary of $I_f$ we must have
$|x| \ll 1$ and $|\xi| \gg 1$). 

\subsection{The regular case}
\label{sec:stab_reg}
Again, we first consider the case $G_1 = \O(1)$ (\ref{def:G1H0}). 
In the slow coordinate $x$ (\ref{fastvar}), i.e. outside the region $I_f$, 
the equation for $u$ reads
\beq
\label{uinv}
(1 - 3U_h^2 - \la + \O(\eps)) u = - U_h v  + \O(\eps^2 u_{xx})
\eeq
(\ref{linstab}), since $V_h(\xi) = \O(\eps)$ on $\mathbb{R}$ (Theorem \ref{th:ex_reg}). 
Thus, $u$ can be expressed in terms of $v$ outside the fast $\xi$-region $I_f$ 
(\ref{def:If}). Using that $U^2_h(\xi;\eps) = 1 + \O(\eps)$ outside $I_f$ 
(Theorem \ref{th:ex_reg}), we find for the $v$-equation of 
(\ref{linstab}) on the slow $x$-scale,
\[
\begin{array}{lcl}
v_{xx} & = &   
\left[2 H(1,0) U_h + \O(\eps)\right] u - \left[H(1,0) + \frac{\pa G}{\pa V}(0) 
- \tau \la  + \O(\eps) \right] v
\\ & = &
\left[\frac{2 H_0}{\la+2} - H_0 - G_1 + \tau \la + \O(\eps)\right] v  +  \O(\eps^2 v_{xx})
\end{array}
\]
(\ref{def:G1H0}). Hence, outside $I_f$
\beq
\label{vxx}
v_{xx} =
\left[\frac{-H_0 \la + \la(\la+2)\tau - G_1(\la+2)}{\la+2} + \O(\eps) \right] v,
\eeq
uniformly in $\xi$.
The $v$-equation is thus at leading order of constant coefficients type.
By (\ref{chareq}) and (\ref{Lambdas}) we have on the $\xi$-scale
\beq
\label{vxixi}
v_{\xi\xi} = 
\left[\frac{Q(\la;0)}{\la+2} + \O(\eps^3) \right]v
= \left[\La^2_{2,3}(\la, \eps) + \O(\eps^3)\right]v,
\eeq 
In order to determine an expression for $t_2(\la,\eps)$, we need to control
the solution $\phi_2(\xi;\la,\eps)$ (Lemma \ref{lem:defphi12t12}) of (\ref{lindyn}). 
This is done in the following Lemma.
\begin{Lem}
\label{lem:t2t3}
For all $\la \in \C_e$ such that $t_1(\la,\eps) \neq 0$
there exist ${\cal O}(1)$ constants $C_{-}, C_{+} > 0$ and a third
meromorphic transmission function $t_3(\la,\eps)$ such that 
\begin{equation}
\phi_2(\xi; \la, \eps) =
\left\{
\begin{array}{rll}
\left[E^-_2(\la) + \O(\eps) \right] e^{\La_2(\la) \xi} + {\cal O}(e^{C_{-}\xi})  
& {\rm for} & \xi < -\frac{1}{\sqrt{\eps}} \\
t_2(\la) E^+_2(\la) e^{\La_2(\la) \xi} + 
t_3(\la) E^+_3(\la) e^{\La_3(\la) \xi} +
{\cal O}(e^{-C_{+}\xi}) & 
{\rm for} & \xi > \frac{1}{\sqrt{\eps}}.
\end{array}
\right.
\label{def:t3}
\end{equation}
Moreover, there exists an ${\cal O}(1)$ constant $C_{f}$  
such that
$||\phi_2(\xi)||\le C_{f}$ for $\xi \in I_f$. The $v$-coordinate of $\phi_2(\xi)$
satisfies $v(\xi)=1+{\cal O}(\sqrt{\eps})$ on $I_f$, so that
\beq
\label{t2+t3}
t_2(\la, \eps) + t_3(\la, \eps) = 1 + \O(\sqrt{\eps}).
\eeq
\end{Lem} 
\begin{proof}
  The behavior of $\phi_2(\xi)$ outside $I_f$ is determined by
  (\ref{vxixi}) and (\ref{uinv}). The approximation (\ref{def:t3}) for
  $\xi < -1/\sqrt{\eps}$ follows from the definition of $\phi_2(\xi)$
  (Lemma \ref{lem:defphi12t12}).  This same Lemma establishes the
  leading order term in (\ref{def:t3}) for $\xi \to \infty$. The
  transmission function $t_3(\la, \eps)$ measures the component of
  $\phi_2(\xi)$ that decays on the slow spatial scale $x$.  Inside
  $I_f$, $v_{\xi\xi} = \O(\eps^2)$ (\ref{linstab}) and
  $\La^2_{2,3}(\la, \eps) = \O(\eps^2)$ (\ref{Lambdas}), so that
  (\ref{t2+t3}) follows.  As in section \ref{sec:Evans} we refrain
  from giving the full analytic details of this result, since these
  are essentially the same as in \cite{gj,dgk2,dgk3}.
\hfill\end{proof}\\
The transmission function $t_2(\la, \eps)$ can be determined by the
methods originally developed in \cite{dgk1}.  We deduce from Lemma
\ref{lem:t2t3} that the total change in $v_{\xi}$ over $I_f$ is given
by \beq
\label{vxi_slow}
\Delta_{\rm slow} v_{\xi} = 2 \eps (t_2(\la) - 1) 
\sqrt{\frac{Q(\la;0)}{\la+2}} + \O(\eps \sqrt{\eps}).
\eeq
This change in $v_{\xi}$ must be an effect of the evolution on the fast $\xi$-scale,
that is given by
\beq
\label{vxi_fast1}
\Delta_{\rm fast} v_{\xi} = \int_{-\frac{1}{\sqrt{\eps}}}^{\frac{1}{\sqrt{\eps}}}
v_{\xi \xi}|_{\{u=u_{\rm in}, v=1\}} d \xi + \O(\eps^2\sqrt{\eps}), 
\eeq
where $u_{\rm in}(\xi)$ is a bounded solution of the inhomogeneous problem
\beq
\label{uineq}
u_{\xi\xi} + (1 - 3U_h^2(\xi;0) - \la) u = - U_h(\xi;0)
\eeq
(recall that $v(\xi) = 1 + \O(\sqrt{\eps})$ in $I_f$).
The transmission function $t_2(\la;\eps)$ is determined by combining
(\ref{vxi_slow}) and (\ref{vxi_fast1}). Since, a priori $\Delta_{\rm slow} v_{\xi} = \O(\eps)$
and $\Delta_{\rm fast} v_{\xi} = \O(\eps \sqrt{\eps})$ we are led to the following
conclusion. 
\begin{Lem}
\label{lem:t2=1}
Consider $\la \in \C_e \cap \{{\rm Re}(\la) > -2 + \delta\}$
for some $\delta >0$ independent of $\eps$.
Let $\la_2^f = -\frac32$ be the second eigenvalue of the limit system (\ref{linfast})
with $V_0 = 0$ (\ref{fastevs}),
and let $\la^+(0)$ and $\la^-(0)$ be the solutions of $Q(\la,0) = 0$ (\ref{chareq}).
Then, $t_2(\la) = 1 + \O(\sqrt{\eps})$ if 
$|\la - \la_2^f|, |\la - \la^+(0)|, |\la - \la^-(0)| = \O(1)$;
$t_2(\la) = 1 + \O(\eps^{\frac12 - \si})$
if $|\la - \la_2^f| = \O(\eps^{\si})$,  
$|\la - \la^+(0)| = \O(\eps^{2 \si})$,  or $|\la - \la^-(0)| = \O(\eps^{2 \si})$
for some $\si \in (0,\frac12)$.  
\end{Lem}

Thus, this Lemma establishes that $t_2(\la, \eps)$ can only be zero in
$\{{\rm Re}(\la)>-2\}$ if $\la \in \C_e$ is $\O(\sqrt{\eps})$ close to
$\la_2^f$ or $\O(\eps)$ close to $\la^+(0)$ or $\la^-(0)$, so we only
have to study $\la$ near these $3$ points to determine the slow-fast
eigenvalues of (\ref{lindyn}).  Note that the fast reduced (scalar)
limit problem has an eigenvalue $\la_2^f = -\frac32$ ((\ref{fastevs}),
since $V_0=V_h(0) \to 0$ as $\eps \to 0$ (Theorem \ref{th:ex_reg})).
We will prove below that $t_2(\la)$ has a (simple) zero close to
$\la_2^f$, i.e. that the fast reduced eigenvalue $\la_2^f$ persists.

However, before going further into the details of the (possible)
existence of eigenvalues near $\la_2^f$, $\la^+(0)$ or $\la^-(0)$, we
formulate a result that is an immediate consequence of Lemma \ref{lem:t2=1}
and that establishes the stability of the wave for
values of $G_1$ and $H_0$ such that the essential spectrum, and hence
$\la^+(0)$ and $\la^-(0)$, is in the negative half-plane and not too
close to the imaginary axis (see Lemma \ref{lem:stabsig_ess}).
\begin{Theorem}
\label{th:reg_stab}
Let $\eps > 0$ be small enough and let $G_1 <0$ and $H_0 + G_1 - 2 \tau < 0$ 
be such that $|G_1|, |H_0 + G_1 - 2 \tau| \gg \eps$.
The spectrum of the eigenvalue problem (\ref{linstab}) associated to the stability of the
solution $(U_h(\xi; \eps), V_h(\xi; \eps))$ consists of a (unique) eigenvalue at $\la = 0$
and a part that is embedded in the region $\{ {\rm Re}(\la) < - \eps \}$. Therefore, 
$(U_h(\xi; \eps), V_h(\xi; \eps))$ is (spectrally) stable.
\end{Theorem}
Note that the operator defined by (\ref{linstab}) is clearly sectorial in this case
(see section \ref{sec:stab_ess}), so that the nonlinear stability of 
$(U_h(\xi; \eps), V_h(\xi; \eps))$ follows by standard arguments (see for instance \cite{henry}).

{\em Proof of Lemma \ref{lem:t2=1}.} 
We first note that indeed $\Delta_{\rm fast} v_{\xi} = \O(\eps \sqrt{\eps})$
and $\Delta_{\rm slow} v_{\xi} = \O(\eps)$, and thus $t_2(\la) = 1 + \O(\sqrt{\eps})$, 
for $\la \in \C_e$ that are not asymptotically close to
the possible degenerations of (\ref{vxi_fast1}) and (\ref{vxi_slow}).

The inhomogeneous function $u_{\rm in}$ may become unbounded
as $\la$ approaches an eigenvalue, $\la_1^f = 0$ or $\la_2^f = -\frac32$,
or the essential spectrum $\si_{\rm ess}^f$ (\ref{siessf}) 
of the linear problem associated to
the fast reduced limit (\ref{linfast}) with $V_0 = 0$. To avoid irrelevant technicalities
near $\si_{\rm ess}^f$ we assume that 
$\la \in \C_e \cap \{{\rm Re}(\la) > -2 + \delta\}$.
The eigenfunction associated to
$\la_1^f$, i.e. $U_{h,\xi}(\xi;0)$ is odd, which implies that the inhomogeneous (and even)
term $U_h(\xi;0)$ satisfies the solvability condition associated to (\ref{uineq})
at $\la = 0$. Hence, $u_{\rm in}$ remains bounded as $\la \to 0$, so that 
$t_2(\la) = 1 + \O(\sqrt{\eps})$ 
also near $\la = 0$ \cite{dgk2,dgk3}.  The eigenfunction associated to
$\la_2^f$ is even, thus $u_{\rm in}$ grows as $1/(\la_2^f - \la)$ as $\la \to \la_2^f$
\cite{titch,dgk2,dgk3}, which implies that $t_2(\la, \eps) - 1 = \O(\eps^{\frac12 -\si})$
if $|\la - \la_2^f| = \O(\eps^{\si})$ for some $\si \in (0, \frac12)$.

The behavior of $t_2(\la)$ near the degenerations of (\ref{vxi_slow}),
i.e. the zeroes $\la^{\pm}(0)$ of $Q(\la;0)$, follows from observing
that $\Delta_{\rm slow} v_{\xi} = (t_2 - 1) \times \O(\eps^{1+\si})$
if $\la$ is $\O(\eps^{2\si})$ close to $\la^{+}(0)$ or to $\la^{-}(0)$
for some $\si \in (0, \frac12)$. 
\hfill $\endproof$

It now follows, by a (standard) winding number
argument \cite{agj,dgk2,dgk3}, that the eigenvalue $\la_2^f$ persists as an eigenvalue of the 
full system (\ref{linstab}) if it is not embedded in the essential spectrum. 
\begin{Lem}
\label{lem:la32} Let $G_1$ and $H_0$ be such that $\sigma_{\rm ess}$ does not intersect
an $\O(\eps^\si)$ neighborhood of $\la_2^f$, for some $\si < \frac12$. Then,
there is an eigenvalue
$\la_2(\eps)$ of (\ref{linstab}) with $\lim_{\eps \to 0} \la_2(\eps) = \la_2^f = -\frac{3}{2}$.
\end{Lem}
\begin{proof} By the assumptions in the Lemma, there exists a contour $K$ in the complex
$\la$-plane that does not intersect $\si_{\rm ess}$,
that encircles an $\O(\eps^{\si})$ neighborhood of $\la_2^f$ 
and that is $\O(\eps^{\si})$ close to $\la_2^f$. 
It follows from Lemma \ref{lem:t2=1} that
$t_2(\la) = 1 + \O(\eps^{\frac12 - \si})$ for $\la \in K$, thus, the winding number
of $t_2(\la)$ over $K$ is $0$. However, $t_2(\la)$ must have a (simple) pole in
the interior of $K$ -- as is observed in the proof of
Lemma \ref{lem:t2=1}. We conclude that $t_2(\la)$ must also have a (simple, real)
eigenvalue in the interior of $K$.
\hfill\end{proof}\\
The possible existence of slow-fast eigenvalues near $\la^+(0)$ or $\la^-(0)$ is much
more subtle. Since such eigenvalues only become relevant to the stability of the
solution $(U_h(\xi;\eps), V_h(\xi;\eps))$ as $G_1$ (or $H_0+G_1-2\tau$) approaches $0$
(Theorem \ref{th:reg_stab}) we will consider this issue in the forthcoming 
sections.   
\begin{Rem} 
\label{rem:eigenfcts}
\rm
The eigenvalues $\la_1(\eps) = 0$ and $\la_2(\eps) \to -\frac32$ as $\eps \to 0$
can be interpreted as `fast' eigenvalues, since they correspond to eigenvalues of
the fast reduced limit problem. However, strictly speaking both eigenvalues also
have the slow-fast structure described in the beginning of this section.

First, we of course know that $\la_1(\eps) = 0$ is an 
eigenvalue -- see also Lemma \ref{lem:lais0}.
Thus it is a zero of $\D(\la, \eps)$. Since $t_2(\la) = 1 + \O(\sqrt{\eps})$ for
$\la$ near $0$, see the proof of Lemma \ref{lem:t2=1}, we conclude that $t_1(0; \eps) \equiv 0$
(note that this in a sense obvious result does not follow directly from Lemma \ref{lem:zeroest1}).
Thus, the solution $\phi_1(\xi; 0, \eps)$ of (\ref{lindyn}) that by construction has 
a purely fast structure
for $\xi \ll -1$,  does not blow up as $e^{\La_1(0,\eps)\xi}$
as $\xi \to \infty$ (Lemma \ref{lem:defphi12t12}).
Nevertheless, the eigenfunction associated to $\la = 0$, $(U_{h,\xi}(\xi), V_{h,\xi}(\xi))$
has a clear slow-fast structure, that it inherits from $(U_{h}(\xi), V_{h}(\xi))$
(Theorem \ref{th:ex_reg}). Hence, $\phi_1(\xi;0, \eps)$ is not the eigenfunction associated to
$\la = 0$. Neither is $\phi_2(\xi;0, \eps)$, since $t_2(0) \neq 0$. It follows that the eigenfunction
associated to $\la = 0$ must be a linear combination of $\phi_1(\xi;0, \eps)$ and 
$\phi_2(\xi;0, \eps)$, and thus that $\phi_1(\xi;0, \eps)$ does not decay as $\xi \to \infty$,
but instead grows linearly (and slowly), as $e^{\La_2(0,\eps)\xi}$ (like $\phi_2(\xi;0, \eps)$).
The linear combination is such that the two growth terms $e^{\La_2(0,\eps)\xi}$ (for $\xi \to \infty$)
cancel.
\\
Second, $\la_2(\eps)$ is not a zero of $t_1(\la)$, although it is asymptotically
close to such a zero, but it is a zero of $t_2(\la)$. Thus, $\phi_2(\xi;\la_2(\eps),\eps)$
is the eigenfunction of (\ref{lindyn}) at $\la = \la_2(\eps)$ (and $\phi_1(\xi;\la_2(\eps),\eps)$
blows up fast, as $e^{\La_1(\la_2(\eps),\eps)\xi}$).
\end{Rem}
 
\subsection{The super-slow case: an example}
\label{sec:stab_suslo_ex}

In the previous section we have seen that the front might destabilize
as $G_1$ approaches $0$ (if we assume that $H_0+G_1 - 2\tau< 0$).  In
this case, Theorem \ref{th:ex_reg} can no longer be used to establish
the existence of the front $(U_h(\xi), V_h(\xi))$. Thus, the question
about the stability of the front is closely related to the
characteristics of the existence problem (as is usual in the analysis
of (traveling) waves, see also \cite{kap98}).  In this section we
consider the bifurcation as $G_1$ approaches $0$. Therefore, we assume
that $H_0- 2\tau< 0$ and $\O(1)$ with respect to $\eps$.  As in
section \ref{sec:exist} we consider in the super-slow case the
simplified system in which the general function $G(V)$ is replaced by
a linear expression: $G(V) = G_1V = -\eps^2 \ga V$ (see Remark
\ref{rem:geng1}).  Note that Theorem \ref{th:reg_stab} a priori
predicts a possible destabilization as $G_1$ becomes $\O(\eps)$, i.e.
already before $G_1 = -\ga \eps^2$, but it will be shown in the next
section that the estimate in Theorem \ref{th:reg_stab} is not sharp,
in the sense that a bifurcation only occurs as $G_1$ decreases to
$\O(\eps^2)$.

One of the main differences between the analysis in this section and
that of the regular case, is the fact that $V_h(\xi)$ is no longer
$\O(\eps)$, i.e. $V_h(\xi)$ does not only contribute to the higher
order terms in the stability analysis of the front solutions.
Nevertheless, we follow the approach of the previous section and
express the solution $u$ of (\ref{linstab}) in terms of $v$, outside
$I_f$ (see (\ref{uinv})), 
\beq
\label{uinv_suslo}
u = -\frac{U_h}{1+V_h-3U_h^2 - \la} v + \O(\eps^2 u_{xx}) =
\left[\frac{U_h}{2(1+V_h) + \la} + \O(\eps^4) \right] v + \O(\eps^2
v_{xx}) 
\eeq 
since $1 + V_h(\xi; \eps) - U_h^2(\xi; \eps) =
\O(\eps^4)$ (see (\ref{u1u2}), recall that
$q^2$ and $G$ are $\O(\eps^2)$ in the super-slow case). This yields
that 
\beq
\label{vxx_suslo_1}
\begin{array}{lcl}
v_{xx} & = & \left\{
2 \left[H(U_h^2,V_h) + \O(\eps^4) \right] U_h u 
- \left[H(U_h^2,V_h) + \O(\eps^4) - \eps^2 \ga - \tau \la \right] v \right\}
\\[2mm] & = &
\left\{ \frac{2 H(U_h^2,V_h) U^2_h}{2(1+V_h) + \la} - H(U_h^2,V_h) 
+ \tau \la + \eps^2\ga + \O(\eps^4)\right\} v + \O(\eps^2 v_{xx})
\\[2mm] & = &
\left\{ \la \left[\tau - \frac{H(1+V_h,V_h)}{2(1+V_h) + \la} \right] 
+ \eps^2\ga + \O(\eps^4)\right\} v + \O(\eps^2 v_{xx}).
\end{array}
\eeq
It follows from section \ref{sec:stab_ess} that one of the `tips' of $\si_{\rm ess}$,
$\la^{+}(0)$, is $\O(\eps^2)$ if $G_1 = \O(\eps^2)$ (and $H_0-2\tau<0$), 
while the other one, $\la^-(0)$, is $\O(1)$ and negative (\ref{la1la2}). Thus,
the destabilization of the front will either be caused by $\si_{\rm ess}$ at
$G_1 = 0 = \ga$, or possibly by a slow-fast eigenvalue $\la$ that is close to $\la^{+}(0)$
(Lemma \ref{lem:t2=1}). Therefore, we introduce $\tla$ by
\beq
\label{deftla}
\la = \eps^2 \tla,
\eeq
which implies that (\ref{vxx_suslo_1}) can also be written as a super-slow system,
\beq
\label{vxx_suslo_2}
v_{xx} = 
\eps^2 \left\{ \tla \left[\tau - \frac{H(1+V_h,V_h)}{2(1+V_h)} \right] 
+ \ga + \O(\eps^2)\right\} v.
\eeq
As in section \ref{sec:ex_suslo}, we first consider the explicit example in which
$H(U^2,V) = H_0 U^2$. Thus, the existence of (several kinds of) front solutions
is established by Theorem \ref{th:ex_suslo}. In this case, the equation for
$v$ is, on the $\xi$-scale, given by
\beq
\label{vxixi_suslo}
v_{\xi\xi} = 
\eps^4 \left[ \tla (\tau - \frac12 H_0) + \ga + \O(\eps^2)\right] v =
\left[\La^2_{2,3}(\la, \eps) + \O(\eps^6)\right]v,
\eeq  
see (\ref{Lambdas}). Note that this equation is of constant coefficients type, and, 
at leading order, the same as in the equation for $v_{\xi\xi}$ in the regular case
(\ref{vxixi}). Hence, we can copy the arguments leading to Lemma \ref{lem:t2t3}
and conclude that the fundamental solution $\phi_2(\xi; \eps^2 \tla, \eps)$ of (\ref{lindyn})
can again be expressed as in (\ref{def:t3}) outside the region $I_f$. Moreover, as in Lemma
\ref{lem:t2t3}, we may conclude that $t_2(\tla, \eps) + t_3(\tla, \eps) = 1 + \O(\sqrt{\eps})$
(\ref{t2+t3}).

We may now proceed as in the preceding section (and as in \cite{dgk2,dgk3}) and determine
$t_2(\tla)$ by measuring the change in the $q=v_{\xi}$-coordinate of $\phi_2(\xi)$ over
the fast field. It follows from (\ref{vxixi_suslo}) that
\beq
\label{vxi_slow_suslo}
\Delta_{\rm slow} v_{\xi} = 2 \eps^2 (t_2(\tla) - 1) 
\sqrt{\tla (\tau - \frac12 H_0) + \ga} + \O(\eps^2 \sqrt{\eps}).
\eeq
Note that we have to assume that $\tla (\tau - \frac12 H_0) + \ga > 0$, i.e.
$\La^2_{2,3}(\la, \eps) > 0$, which is a natural assumption, since 
\beq
\label{tipsiess}
\tla_{\rm tip} = \tla^+(0) = -\frac{2\ga}{2\tau - H_0} < 0
\eeq 
determines the `tip' of $\si_{\rm ess}$ (recall that $H_0 - 2\tau < 0$),
i.e. $t_2(\tla)$ is not defined if $\tla \leq \tla_{\rm tip}$. By definition, 
$\Delta_{\rm fast} v_{\xi}$ is given by (\ref{vxi_fast1}). 
Since, at leading order $V_h(\xi) = V_h(0) = v_0$ and
$U_h(\xi) = u_0(\xi; v_0)$ (uniformly) in $I_f$(Theorem \ref{th:ex_suslo}), 
and since $u_0(\xi; v_0)$
decays exponentially fast on the (fast) $\xi$-scale, it follows that
\beq
\label{vxi_fast_suslo}
\Delta_{\rm fast} v_{\xi} = \eps^2 H_0 \int^{\infty}_{-\infty}
\left\{ 2\left[2u_0^2(\xi;v_0) - 1 - v_0 \right] u_0(\xi;v_0)u_{\rm in}(\xi;v_0) 
- u_0^2(\xi;v_0) \right\} d\xi
+ \O(\eps^2 \sqrt{\eps}), 
\eeq 
where $u_{\rm in}(\xi;v_0)$ is the (uniquely determined) bounded solution of
the inhomogeneous problem
\[
u_{\xi\xi} + (1 + v_0 - 3u_0^2(\xi;v_0)) u = - u_0(\xi;v_0).
\]
Since we already know one solution of the homogeneous problem, $u(\xi) = u_{0,\xi}(\xi;v_0)$,
we can determine $u_{\rm in}(\xi;v_0)$ explicitly,
\beq
\label{uin}
u_{\rm in}(\xi;v_0) = \frac{1}{2(1+v_0)} \left(u_0(\xi;v_0) + \xi u_{0, \xi}(\xi;v_0) \right).
\eeq
Thus, by (\ref{u0p0}), $\Delta_{\rm fast} v_{\xi}$ can be computed explicitly (at leading order),
\[
\Delta_{\rm fast} v_{\xi} = - \eps^2 H_0 \sqrt{2} \sqrt{1+v_0} + \O(\eps^2 \sqrt{\eps}).
\]
Combining this with (\ref{vxi_slow_suslo}) yields an explicit expression for $t_2(\tla)$,
\beq
\label{t2_ex}
t_2(\tla, \eps) = 1 - H_0 \sqrt{ \frac{2(1+v_0)}{\tla (\tau - \frac12 H_0) + \ga} }
+ \O(\sqrt{\eps})
\eeq
for $\tla > \tla_{\rm tip}$ (\ref{tipsiess}). It follows that $t_2(\tla) \geq 1 + \O(\sqrt{\eps})$ 
for $H_0 \leq 0$ and $t_2(\tla) < 1 + \O(\sqrt{\eps})$ for $H_0 > 0$. Hence, $t_2(\tla)$ cannot
have zeroes if $H_0 \leq 0$. In other words, there cannot be an eigenvalue near the tip
of the essential spectrum in case {\bf (ii)} of Theorem \ref{th:ex_suslo}. 
On the other hand, $t_2(\tla)$ can be $0$ for
$H_0 >0$, i.e. in case {\bf (i)} of Theorem \ref{th:ex_suslo} 
there indeed is a `new' slow-fast 
eigenvalue of (\ref{linstab}), it is given by
\beq
\label{laedge_ex}
\la_{\rm edge} = \eps^2 \tla_{\rm edge} = \frac{-2 \ga +
  H_0^2(1+v_0)}{2\tau - H_0} \eps^2 + \O(\eps^2\sqrt{\eps}) > \eps^2
\tla_{\rm tip} = \la_{\rm tip} \eeq (\ref{tipsiess}). Note that the
eigenvalue $\la_{\rm edge}$ merges with $\la_{\rm tip}$, and thus with
$\si_{\rm ess}$ as $H_0 \downarrow 0$. This is of course a leading
order result, the accuracy of our analysis only allows us to conclude
that $|\la_{\rm tip} - \la_{\rm edge}| \leq \O(\eps^2\sqrt{\eps})$ as
$H_0 \downarrow 0$, and that $\la_{\rm edge}$ does not exist for $H_0
< 0$. Nevertheless, we conclude that $\la_{\rm edge}$ appears from the
essential spectrum as $H_0$ increases through $0$. In other words,
$\la_{\rm edge}$ is created, or annihilated, by an edge bifurcation.
Note that the new eigenvalue appears exactly as $\si_{\rm ess}$ becomes
complex valued (section \ref{sec:stab_ess}, Figure \ref{fig:essspec}).

The existence, or non-existence of $\la_{\rm edge}$ is crucial to the
character of the destabilization (see also the numerical simulations
in section \ref{sec:simdis}).  For $H_0 < 0$, the front solution
$(U_h(\xi), V_h(\xi))$ destabilizes as $\ga$, or equivalently $G_1$,
crosses through $0$. The destabilization is due to the essential
spectrum, which implies that also the `background states' $(U(x,t),
V(x,t)) \equiv (\pm 1, 0)$ destabilize at $\ga = 0$.  However, in the
case $H_0 > 0$ the eigenvalue is $\la_{\rm edge}$ is $\eps^2
H_0^2(1+v_0)/(2\tau - H_0)$ ahead of $\si_{\rm ess}$
(\ref{laedge_ex}), in the sense that it reaches the axis Re$(\la) =
0$ before $\si_{\rm ess}$ as $\ga > 0$ decreases to $0$.  Thus, if
$H_0 > 0$ the front solution $(U_h(\xi), V_h(\xi))$ destabilizes by an
element of the discrete spectrum of (\ref{linstab}) at $\ga_{\rm
  double} = \frac12 H_0^2 (1+v_0) + \O(\sqrt{\eps}) > 0$.  As a
consequence, the background states $(\pm 1, 0)$ remain stable as
$(U_h(\xi), V_h(\xi))$ destabilizes for $H_0 > 0$, contrary to the
case $H_0 < 0$.  The bifurcation at $\ga_{\rm double}$ is associated
to the saddle-node bifurcation of heteroclinic orbits described in
Theorem \ref{th:ex_suslo}.
\begin{Theorem}
\label{th:edge_SN_ex}
Assume that $G(V) = -\eps^2 \ga$, $H(U^2,V) = H_0 U^2$, $H_0 - 2\tau < 0$ and $\O(1)$, and 
that $\eps > 0$ is small enough.
\\
{\bf (i)} Let $(U_h^{+,1}(\xi), V_h^{+,1}(\xi))$ and $(U_h^{+,2}(\xi), V_h^{+,2}(\xi))$ 
be the two types of heteroclinic
front solutions that exist for $H_0 > 0$ and $\ga \geq \ga_{\rm double} = \frac32 H_0^2 + \O(\sqrt{\eps})$ 
with, at leading order, 
$0 < V_h^{+,1}(0) = v_1 \leq 2 \leq v_2 = V_h^{+,2}(0)$ (Theorem \ref{th:ex_suslo}). 
The front solution $(U_h^{+,1}(\xi), V_h^{+,1}(\xi))$ is asymptotically stable for 
$\ga > \ga_{\rm double}$,
the front $(U_h^{+,2}(\xi), V_h^{+,2}(\xi))$ unstable; $(U_h^{+,1}(\xi), V_h^{+,1}(\xi))$ destabilizes 
by an element of the discrete spectrum, $\la_{\rm edge}$, at 
$\ga = \ga_{\rm double}$ and merges with $(U_h^{+,2}(\xi), V_h^{+,2}(\xi))$ 
in a saddle-node bifurcation of heteroclinic orbits.
\\ 
{\bf (ii)} Let $(U_h^+(\xi), V_h^+(\xi))$ be a heteroclinic
front solution that exist for $H_0 < 0$ and (all) $\ga > 0$ (Theorem \ref{th:ex_suslo});
$(U_h^+(\xi), V_h^+(\xi))$ is asymptotically stable for all $\ga > 0$, it is destabilized at
$\ga = 0$ by the essential spectrum $\si_{\rm ess}$.   
\end{Theorem} 
\begin{Rem} \rm
As in the regular case, spectral stability implies asymptotic nonlinear stability in
this super-slow case, since the linear operator associated to the stability problem remains sectorial
as long as $\eps > 0$.  
\end{Rem}

{\em Proof of Theorem \ref{th:edge_SN_ex}.} We first note that the condition $H_0 - 2\tau < 0$ and $\O(1)$
determines that $\si_{\rm ess}$ can only cross, or come close to, the Re$(\la) = 0$-axis 
at $\la = 0$ (Lemma \ref{lem:stabsig_ess} with $G_1 = \O(\eps^2)$). 
\\ 
{\bf (i)} The eigenvalue `in front of' $\si_{\rm ess}$, $\la^{1,2}_{\rm edge}(v_{1,2})$, 
is given by (\ref{laedge_ex}),
where $v_0>0$ is a solution of $9 \ga v^2 = 2 H_0^2 (1+v)^3$, and
$v_0 = v_1 \leq 2$ (at leading order) for 
$(U_h^{+,1}(\xi), V_h^{+,1}(\xi))$, while $v_0 = v_2 \geq 2$ (at leading order) for 
$(U_h^{+,2}(\xi), V_h^{+,2}(\xi))$ -- Theorem \ref{th:ex_suslo}. Thus, by (\ref{laedge_ex})
|$\la^{1}_{\rm edge}(v_{1}) < 0$
and $\la^{2}_{\rm edge}(v_{2}) > 0$ if
$\ga < \ga_{\rm double} = \frac32 H_0^2 + \O(\sqrt{\eps})$. As a consequence, 
$\la^{1}_{\rm edge}(v_{1}) \uparrow 0$
and $\la^{2}_{\rm edge}(v_{2}) \downarrow 0$ as $\ga \downarrow \ga_{\rm double}$, at which
the saddle-node bifurcation takes place.. 
\\
{\bf (ii)} We have already shown that there can be no eigenvalues in front of the tip of $\si_{\rm ess}$.
Therefore, the statement of the Theorem follows.  \hfill $\endproof$

\begin{Rem}
\label{rem:edgeeven} \rm
Since $t_2(\la) = 0$, the slow-fast eigenfunction associated to the bifurcation at $\ga = \ga_{\rm double}$
is given by $\phi_2(\xi)$.
It follows from Lemmas \ref{lem:evenodd} and \ref{lem:lais0} that the $u$-component of $\phi_2$
is odd, and the $v$-component even, as function of $\xi$. 
\end{Rem}

\subsection{Bifurcations in the general super-slow problem}
\label{sec:stab_suslo_bif}

We now consider the stability of a front solution in the general
super-slow limit. Thus, we assume we have established the existence of
a front $(U_h(\xi), V_h(\xi))$ for a certain given function $H(U^2,V)$
(Theorem \ref{th:gen_suslo}). To analyze its stability, we again
try to determine $t_2(\la)$ by measuring $\Delta_{\rm fast} v_{\xi}$
and $\Delta_{\rm slow} v_{\xi}$.

In order to determine $\Delta_{\rm slow} v_{\xi}$ we follow the
derivation of (\ref{vxx_suslo_2}) in the previous section. Hence, we
again conclude that non-trivial eigenvalues near $0$ are only possible
for $\la = \O(\eps^2)$, thus we again introduce $\tla$ (\ref{deftla})
(see also the proof of Theorem \ref{th:destab} for more details on the
necessity of this scaling).  Note that both $G_1$ and $\la$ are now
$\O(\eps^2)$, thus, we can immediately obtain a leading order
expression for $\Delta_{\rm fast} v_{\xi}$ in terms of $H(U^2,V)$,
\beq
\label{vxi_fast_gen}
\begin{array}{ccc}
\Delta_{\rm fast} v_{\xi} & = &
\eps^2 \int_{-\infty}^{\infty} \left\{ 
2 \left[H(u_0^2,v_0) - (1+v_0-u_0^2)
\frac{\pa H}{\pa U^2}(u_0^2,v_0)\right] u_0 u_{\rm in} \right. \\
& & \left. \; \; \; \; \; \; \; \; - \left[H(u_0^2,v_0) - (1+v_0-u_0^2)
\frac{\pa H}{\pa V}(u_0^2,v_0) \right] \right\} d\xi + \O(\eps^2 \sqrt{\eps})
\end{array}
\eeq
(\ref{linstab}), where 
$u_{\rm in}(\xi)$ is given in (\ref{uin}) -- recall that 
$v = 1 + \O(\sqrt{\eps})$ in $I_f$. 
As in the previous section, we have approximated $U_h(\xi)$ by
$u_0(\xi;v_0)$ (\ref{u0p0}), $V_h(\xi)$ by $v_0$ and $I_f$ by $\mathbb{R}$
(Theorem \ref{th:gen_suslo}).
Note that the integral converges and that $\Delta_{\rm fast} v_{\xi}$ is (at leading order)
independent of $\ga$ and $\tla$.

It is in principle possible to determine $\Delta_{\rm slow} v_{\xi}$
in terms of $t_2(\la)$
from (\ref{vxx_suslo_2}), however, this equation is in general not of constant coefficients
type (unlike for the example problem in section \ref{sec:stab_suslo_ex}). If we introduce the
super-slow coordinate $X$ by $X = \eps x = \eps^2 \xi$, we can write (\ref{vxx_suslo_2}) as
\beq
\label{vXX}
v_{XX} = \left\{ \tla \left[\tau -
    \frac{H(1+V_h(X),V_h(X))}{2(1+V_h(X))} \right] + \ga +
  \O(\eps^2)\right\} v, \eeq i.e. the functions $V_h(X)$ introduce
explicit $X$-dependent terms in the equation (section
\ref{sec:gen_suslo}, $V_h(X)$ behaves as $e^{\mp \sqrt{\ga} X}$ on
$\M^{\pm}_{\eps}$). Nevertheless, we can in principle determine the
$v$-components of the solution $\phi_2(\xi)$ of (\ref{lindyn}) outside
the fast region $I_f$. However, the analysis is much less transparent.
For instance, the decomposition (\ref{def:t3}) as in Lemma
\ref{lem:t2t3} now only holds for $X \gg 1$, therefore the relation
between $t_3(\tla)$ and $t_2(\tla)$ that is obtained from the value of
$v$ in $I_f$ will in general be more complicated than in
(\ref{t2+t3}). Moreover, $\tla \left[\tau -
  \frac{H(1+V_h(X),V_h(X))}{2(1+V_h(X))} \right] + \ga$ might change
sign as function of $X$, so that the solution $v(X)$ of (\ref{vXX})
can have oscillatory parts.

Thus, we conclude that it is not a straightforward extension of the
approach in previous section to determine $t_2(\tla)$ for general
values of $\tla$. It should also be noted that a similar problem
occurs in the regular case, in the study of possible eigenvalues near
$\la^{\pm}(0)$ (Lemma \ref{lem:t2=1}). If one introduces $\tla^{\pm}$
by $\la = \la^{\pm}(0) + \eps \tla^{\pm}$, and derives the leading
order equation for $v_{xx}$ (\ref{vxx}) in this case, then one finds
an equation like (\ref{vXX}), i.e. an equation with spatially
dependent coefficients (these $x$-dependent terms originate from the
$\O(\eps)$ corrections corresponding to $V_h(x) = \O(\eps)$ in
(\ref{uinv}) and (\ref{vxx})).
Hence, at this point it is not yet possible to
determine in full detail whether or not eigenvalues exist near the
tips of $\si_{\rm ess}$ for general nonlinearities $H(U^2,V)$ and
general $\la$. Moreover, it is also not possible to explicitly describe
how and when eigenvalues appear from, or disappear into,
$\si_{\rm ess}$. On the other hand, it is clear from (\ref{vxi_fast_gen})
and (\ref{vXX}) that the number of zeros of $t_2(\tla)$ depends
(for instance) on $H_0$. It thus follows that eigenvalues will be 
created/annihilated near the tip of $\si_{\rm ess}$ in the general case
(like in the example system considered in the previous section).
The analysis of eigenvalues near the tip of $\si_{\rm ess}$ is therefore
a continuing subject of research (in progress, see also
section \ref{sec:simdis}).

Nevertheless, the value $\la = \tla = 0$ is, of
course, especially relevant for the stability analysis of the front,
and, equation (\ref{vXX}) is again of constant coefficients type at
leading order for this special value of $\la$. Hence, for $\la = 0$ we
can obtain the equivalent of Lemma \ref{lem:t2t3}, so that it follows
that 
\beq
\label{vxi_slow_la=0}
\Delta_{\rm slow} v_{\xi}|_{\la = 0} = 2 \eps^2 (t_2(0) - 1) 
\sqrt{\ga} + \O(\eps^2 \sqrt{\eps}).
\eeq
Note that eventually it becomes clear at this point why the choice 
$G_1 = -\eps^2 \ga$ is the most relevant scaling of $G_1$. With this scaling
the `jumps' $\Delta_{\rm slow} v_{\xi}$ and $\Delta_{\rm fast} v_{\xi}$ (\ref{vxi_fast_gen})
are of the same magnitude in $\eps$ at $\la = 0$. Therefore, $t_2(0, \eps)$ is asymptotically
close to $1$ for all $G_1$ with $|G_1| \gg \eps^2$ -- see Lemma \ref{lem:t2=1} and its proof.
Thus, the stability problem (\ref{linstab}) can only have a double eigenvalue at $0$ 
if $G_1 = \O(\eps^2)$. 
This establishes a significant link between the stability analysis and the existence
analysis of section \ref{sec:exist}, since it is clear from the analysis there that
the scaling $G_1 = \O(\eps^2)$ is also the most relevant scaling for the (super-slow)
existence problem (Remark \ref{rem:sigdeg}). Moreover, this link is even much more explicit.
\begin{Theorem}
\label{th:t2=0SN}
Assume that $G(V) = -\eps^2 \ga$, that $H_0 - 2\tau < 0$ and $\O(1)$,
and that $\eps > 0$ is small enough. Let the front
solution $(U_h(\xi;\eps),V_h(\xi;\eps))$ be a heteroclinic solution 
that corresponds to an intersection $T^-_o \cap W^u(-1,0,0,0)|_{\M_\eps^-}$ as described in
Theorem \ref{th:gen_suslo}.  
The stability problem associated to the front solution has a double eigenvalue
at $\la = 0$ if and only if the intersection $T^-_o \cap W^u(-1,0,0,0)|_{\M^-_\eps}$ 
is non-transversal. If the intersection $T^-_o \cap W^u(-1,0,0,0)|_{\M^-_\eps}$ is a second order
contact, then the front bifurcates at
\beq
\label{gaedge_gen}
\begin{array}{rl}
0 < \ga_{\rm double} = \frac{1}{4(1+v_0)^2} & \left[  
\int_{-\infty}^{\infty} (1+v_0-u_0^2)H(u_0^2,v_0) d\xi \right. \\
&
\left.
+ 2 \int_{-\infty}^{\infty} (1+v_0-u_0^2)[u_0^2 \frac{\pa H}{\pa U^2}(u_0^2,v_0) +
(1+v_0) \frac{\pa H}{\pa V}(u_0^2,v_0)] d\xi \right]^2
\end{array}
\eeq
by merging with another front
solution in a saddle-node bifurcation of heteroclinic orbits. 
\end{Theorem}
\begin{proof} First, we recall from section \ref{sec:gen_suslo}
that a heteroclinic connection
that corresponds to the intersection of $W^u(-1,0,0,0)|_{\M_\eps^-} = 
\{q = \eps \sqrt{\ga} v\}$ and $T^-_o$ is determined by (\ref{defhets_gen}).
This is of course a leading order approximation. 
In the proof of this Theorem we refrain from mentioning this obvious fact at several places.
To determine the $v_0$-dependence of the right hand side of this relation, we define $w_0(\xi)$
as the (monotonically increasing) heteroclinic solution of $\ddot w + (1-w^2) w = 0$. 
It follows that
\beq
\label{defw0}
u_0(\xi;v_0) = \sqrt{1+v_0} w_0(\sqrt{1+v_0} \xi), \; \; w_0(t) = \tanh \sqrt{\frac12} t
\eeq
(\ref{u0p0}). Replacing $u_0(\xi;v_0)$ by $w_0(t)$ in (\ref{defhets_gen}) yields
\beq
\label{hets_gen(v0)}
\sqrt{\ga} v_0 = \frac12 \sqrt{1+v_0} \int_{-\infty}^{\infty}
(1-w_0^2)H((1+v_0)w_0^2,v_0) dt.
\eeq
Thus, $T^-_o \cap W^u(-1,0,0,0)|_{\M_\eps^-}$ is non-transversal if (\ref{defhets_gen}) holds
and
\beq
\label{nontrans}
\begin{array}{rl}
\sqrt{\ga} = &
\frac12 \frac{\pa}{\pa v_0} \left\{ \sqrt{1+v_0} \int_{-\infty}^{\infty}
(1-w_0^2)H((1+v_0)w_0^2,v_0) dt \right\}
\\
= &
\frac{1}{4 \sqrt{1+v_0}} 
\int_{-\infty}^{\infty} (1-w_0^2)H((1+v_0)w_0^2,v_0) dt 
\\
&
+ \frac12 \sqrt{1+v_0} 
\int_{-\infty}^{\infty} (1-w_0^2)[w_0^2 \frac{\pa H}{\pa U^2}((1+v_0)w_0^2,v_0) +
\frac{\pa H}{\pa V}((1+v_0)w_0^2,v_0)] dt
\\
= &
\frac{1}{2(1+v_0)}  
\int_{-\infty}^{\infty} (1+v_0-u_0^2)H(u_0^2,v_0) d\xi 
\\
&
\frac{1}{1+v_0}
\int_{-\infty}^{\infty} (1+v_0-u_0^2)[u_0^2 \frac{\pa H}{\pa U^2}(u_0^2,v_0) +
(1+v_0) \frac{\pa H}{\pa V}(u_0^2,v_0)] d\xi,
\end{array}
\eeq
by re-introducing $u_0(\xi;v_0)$. Note that (\ref{gaedge_gen}) follows from this equation.
The expression for $t_2(0,\eps)$ is determined by (\ref{vxi_fast_gen}), (\ref{vxi_slow_la=0})
and (\ref{uin}),
\[
t_2(0,\eps) = 1 - \frac{{\cal I}_1 + {\cal I}_2 + {\cal I}_3}{2 \sqrt{\ga}(1+v_0)} + \O(\sqrt{\eps}),
\]
where
\beq
\label{I123}
\begin{array}{rcl}
{\cal I}_1 & = &
\int_{-\infty}^{\infty} (1+v_0-u_0^2)H(u_0^2,v_0) d\xi,
\\ 
{\cal I}_2 & = &
\int_{-\infty}^{\infty} (1+v_0-u_0^2)[u_0^2 \frac{\pa H}{\pa U^2}(u_0^2,v_0) +
(1+v_0) \frac{\pa H}{\pa V}(u_0^2,v_0)] d\xi,
\\ 
{\cal I}_3 & = & 
\int_{-\infty}^{\infty} [(1+v_0-u_0^2)\frac{\pa H}{\pa U^2}(u_0^2,v_0) - H(u_0^2,v_0)]
\xi u_0 u_{0,\xi} d\xi.
\end{array}
\eeq
We find by partial integration that
\[
{\cal I}_3 = \int_{-\infty}^{\infty} \frac12 \xi \frac{\pa}{\pa \xi}
[(1+v_0-u_0^2) H(u_0^2,v_0)] d \xi = -\frac12 {\cal I}_1,
\]
which implies that
\[
t_2(0,\eps) = 1 - \frac{{\cal I}_1 + 2{\cal I}_2}{4 \sqrt{\ga}(1+v_0)} 
+ \O(\sqrt{\eps}),
\]
so that we can conclude by (\ref{I123}) that $t_2(0, \eps) = 0$ is equivalent 
to the non-transversality
condition (\ref{nontrans}). Hence, a double eigenvalue of (\ref{linstab}) coincides with a
saddle-node bifurcation of heteroclinic orbits, unless the tangency between
$T^-_o$ and $W^u(-1,0,0,0)|_{\M_\eps^-}$ is degenerate.
\hfill\end{proof}\\
Finally, we can turn to the question about the character of the 
destabilization of the regular front
solution, that has been studied in sections \ref{sec:ex_reg} and
\ref{sec:stab_reg}, as $G_1$ approaches $0$.
In order to do so, we first note that the existence problem for the regular case can be recovered 
from that of the singular limit by re-introducing $G_1 = -\ga \eps^2$ 
in the existence condition (\ref{defhets_gen}). This implies that $v_0$ must become
$\O(\eps)$ and that
\beq
\label{limreg}
\sqrt{-G_1} v_0 = \eps \frac12 \int_{-\infty}^{\infty}
(1-u_0^2)H(u_0^2,0) d\xi + \O(\eps \sqrt{\eps}), 
\eeq 
which is
equivalent to (\ref{appr_vheps}) in Theorem \ref{th:ex_reg}. Thus, the
structure of the front $(U_h(\xi),V_h(\xi))$ as function of $G_1
\uparrow 0$ can be determined by tracing the intersection $T^-_o \cap
W^u(-1,0,0,0)|_{\M_\eps^-}$ in the super-slow limit as
$W^u(-1,0,0,0)|_{\M_\eps^-} = \{q = \eps \sqrt{\ga} v \}$ goes down from
being almost vertical ($G_1 = \O(1)$, $\ga = \O(1/\eps^2)$) to
horizontal ($G_1 = \ga = 0$). Note that this process determines a
unique `regular' element in the intersection $T^-_o \cap
W^u(-1,0,0,0)|_{\M_\eps^-}$, all other elements of $T^-_o \cap
W^u(-1,0,0,0)|_{\M_\eps^-}$ do not persist in the regular limit $\ga =
\O(1/\eps^2)$ (here, we do not pay attention to possible heteroclinic
connections that have $v_0 \gg 1$ as $\ga \gg 1$).  It depends on the
sign of $\frac12 \int_{-\infty}^{\infty} (1-u_0^2)H(u_0^2,0) d\xi$
whether $v_0$ will be positive or negative (\ref{limreg}), i.e.
whether the regular intersection $T^-_o \cap W^u(-1,0,0,0)|_{\M_\eps^-}$
travels through the first or through the third quadrant of the
$(v,q)$-plane as $\ga$ decreases. Since $H(U^2,V)$ is smooth, we can
make a distinction between two different types of behavior:


\begin{description}
\item{Type D} The regular element of $T^-_o \cap
  W^u(-1,0,0,0)|_{\M_\eps^-}$ merges at a certain critical value of
  $G_1 = - \eps^2 \ga < 0$ with another element of $T^-_o \cap
  W^u(-1,0,0,0)|_{\M_\eps^-}$ in a saddle-node bifurcation of
  heteroclinic orbits.
\item{Type E} The regular element of $T^-_o \cap
  W^u(-1,0,0,0)|_{\M_\eps^-}$ exists up to the limit $G_1 = 0$.
\end{description}
Note that $T^-_o$ approaches $(-1,0)$
as $v_0 \downarrow -1$ (\ref{hets_gen(v0)}), so that an element of 
$T^-_o \cap W^u(-1,0,0,0)|_{\M_\eps^-}$ can
only reach the singular region $\{v_0 \leq -1\}$ at $\ga = 0$, which 
indeed implies that there can
only be orbits of type D and E in the third quadrant. We can now 
describe the destabilization of
the regular fronts as $G_1$ approaches $0$. 
\begin{Theorem}
\label{th:destab}
Assume that $G(V) = -\eps^2 \ga$, that $H_0 - 2\tau < 0$ and $\O(1)$,
and that $\eps > 0$ is small enough. 
Consider the heteroclinic front solution $(U_h(\xi),V_h(\xi))$ 
determined in Theorem \ref{th:ex_reg}
for $G_1 < 0$ and $\O(1)$ and in Theorem \ref{th:gen_suslo} for $G_1 = \O(\eps^2)$.
If the front is of type D as $G_1$ becomes $\O(\eps^2)$, then
it is asymptotically stable up to $G_1 = -\eps^2 \ga_{\rm double} < 0$ (\ref{gaedge_gen})
and it is destabilized by a (discrete) eigenvalue
through a saddle-node bifurcation of heteroclinic
orbits. A front solution of type E is stable up to $G_1 = 0$ and it is destabilized 
by the essential spectrum.
\end{Theorem}

Thus, the destabilization of a regular front solution in the limit 
$G_1 \uparrow 0$ is completely determined
by the geometrical structure of $T^-_o \cap W^u(-1,0,0,0)|_{\M_\eps^-}$ in the super-slow limit. 
Note that Figure \ref{fig:gensuslo} presents examples of type D and type E behavior.    
\begin{proof}
The proof of this Theorem is a bit more subtle than a priori might be expected, 
since in general we do not have control over
the eigenvalues of (\ref{linstab}) near the tip of $\si_{\rm ess}$ 
(see also Remark \ref{rem:manyev}),
except that these eigenvalues must be $\O(\eps^2)$ close to $\si_{\rm ess}$ (see also below). 
Thus, for instance the following scenario for a type D orbit
might be possible as $\ga$ decreases to $\ga_{\rm double}$: 
two eigenvalues bifurcate (subsequently) from $\si_{\rm ess}$ (as real eigenvalues), merge and
become a pair of complex eigenvalues. This pair crosses through the Re$(\la) = 0$ axis
at $\ga_{\rm Hopf} > \ga_{\rm double}$, and touch down again on the real axis. At 
$\ga_{\rm double}$ one of these eigenvalues returns to Re$(\la) = 0$. Thus, in this scenario,
there already exists an unstable eigenvalue at $\ga = \ga_{\rm double}$; moreover, the front
destabilizes by a Hopf bifurcation at $\ga_{\rm Hopf} > \ga_{\rm double}$.

Let us first note that a destabilization by a Hopf bifurcation is the only alternative
to the statements of the Theorem, since eigenvalues either move through $0$, or (in pairs) through
the Re$(\la) = 0$ axis. If we can show that a Hopf bifurcation cannot occur for
$\ga > \ga_{\rm double}$, then it is clear that for type D orbits
$\la_{\rm edge} < 0$ for $\ga > \ga_{\rm double}$ 
and that there is no unstable spectrum at $\ga = \ga_{\rm double}$
(this follows from Theorem \ref{th:reg_stab}: if $\ga$ is $\gg \O(1/\eps)$,
all non-trivial eigenvalues must be in $\{$Re$(\la) < -\eps\}$, hence, by decreasing
$\ga$, there is one eigenvalue, $\la_{\rm edge}$, that is the first to reach $0$;
this happens at the saddle-node bifurcation (Theorem \ref{th:t2=0SN}), i.e. at 
$\ga = \ga_{\rm double}$). Thus, the front is stable for $\ga > \ga_{\rm double}$.
The same argument can be used to establish the non-existence of unstable spectrum
for type E orbits, if there are no Hopf bifurcations possible.

To show that there cannot be Hopf bifurcations (for $H_0 - 2\tau < 0$ and $\O(1)$, see
section \ref{sec:simdis}), we
first ascertain that $\la$ must be $\O(\eps^2)$, i.e. that (\ref{deftla})
is the correct scaling. This follows by the same arguments as in the proof of Lemma
\ref{lem:t2=1}. If $\Delta_{\rm slow} v_{\xi} \gg \Delta_{\rm fast} v_{\xi}$, then
there cannot be an eigenvalue. Thus, it follows from (\ref{vxx_suslo_1}) that
$|\la|$ must indeed be $\O(\eps^2)$ near $\la^+(0)$. Hence, even if there is a Hopf
bifurcation, it will be $\O(\eps^2)$ close to $0$. Next, we realize that this situation 
is covered by (\ref{vxi_fast_suslo}) for the jump through the fast
field, thus, $\Delta_{\rm fast} v_{\xi}$ is real
(at leading order), independent of $\tla$. 
However, it follows from (\ref{vXX}) that $\Delta_{\rm slow} v_{\xi}$
cannot be real if $\tla$ is complex-valued. Hence, there cannot be a Hopf bifurcation
$\O(\eps^2)$ close to $\la = 0$.
\hfill\end{proof} 
\begin{Rem}
\label{rem:manyev}
\rm
By the same geometrical arguments (that are based on Theorem \ref{th:t2=0SN}) we can
describe the character of the bifurcations as function $\ga$ in the stability problem associated
to a heteroclinic orbit that corresponds to a 
non-regular element of 
$T^-_o \cap W^u(-1,0,0,0)|_{\M_\eps^-}$.
However, it should be noted that, in general, we do not have enough information on
the spectrum of (\ref{linstab}) to establish
the stability of such a front, since we did not determine all possible eigenvalues.
In general, we cannot exclude the possibility that various eigenvalues have
bifurcated from the essential spectrum for these fronts (in fact, the possible
oscillatory character of a solution $v(X)$ of (\ref{vXX}) strongly suggests that
this can happen). Nevertheless, we may for instance conclude that if the regular orbit
is of type D, then it merges with a non-regular orbit at $\ga_{\rm double}$ that
is unstable for any $\ga > \ga_{\rm double}$ for which it exists. 
\end{Rem}
\begin{Rem} \rm
  The most simple example one can consider is $H(U^2,V) \equiv H_0$.
  This corresponds to the case in which the function $F(U^2,V)$ in
  (\ref{most}) is (the most general) linear function of $U^2$ and $V$
  with parameters $G_1$ and $H_0$ (i.e. $F(U^2,V) = H_0 + (H_0+G_1)V
  -H_0U^2$, recall that $F(1,0)$ must be $0$). In this case, $T^-_o$
  is given by $\{q = 2 \eps H_0 \sqrt{1+v_0}+\O(\eps^2) \}$, so that
  $W^u(-1,0,0,0)|_{\M_\eps^-}$ can never be tangent to $T^-_o$. Hence in
  this case, there is a uniquely determined front solution of type E
  for any $H_0 \neq 0$ and $G_1 < 0$, i.e. the front solution is
  stable up to $G_1 = 0$ and it is destabilized by the essential
  spectrum.
\end{Rem}
\begin{Rem}
\label{rem:degH} \rm
  We did not consider the degenerate case in which $H(U^2,V)$ is such
  that $H(1+V,V) \equiv 0$ (section \ref{sec:intro}), i.e. functions $H$ such that
  $H(U^2,V) = (1+V-U^2)\tilde{H}(U^2,V)$ for some smooth function
  $\tilde{H}$. In a sense, this is a much more simple problem, for
  instance since in the super-slow limit, the stability problem in the
  slow field is automatically of constant coefficients type (at
  leading order), see (\ref{vxx_suslo_1}), (\ref{vXX}). Moreover, it
  is also clear from these same relations that we can find $\O(1)$
  instead of $\O(\eps^2)$ eigenvalues in this case if $\tau =
  \O(\eps^2)$. In fact, the situation is very much like the stability
  analysis of (homoclinic) pulses in mono-stable systems in
  \cite{dgk2,dgk3}. For instance, as in \cite{dgk2,dgk3}, potential
  eigenvalues are no longer `slaved' to the tips of the essential
  spectrum or to the eigenvalues of the fast reduced limit (Lemma
  \ref{lem:t2=1}). Moreover, the `natural' persistence result of Lemma
  \ref{lem:la32} is also not valid in this case, in general.
\end{Rem}
\section{Simulations and Discussion}
\label{sec:simdis}
\setcounter{Theorem}{0} \setcounter{equation}{0}
\subsection{Simulations}
\begin{figure}[ht]\label{desblowup}
\begin{center}
\psfrag{x}{$\xi$}
\psfrag{t}{$t$}
\psfrag{u}{$U$}
\psfrag{v}{$V$}
\subfigure[$U$-component]{\includegraphics[scale=.4]{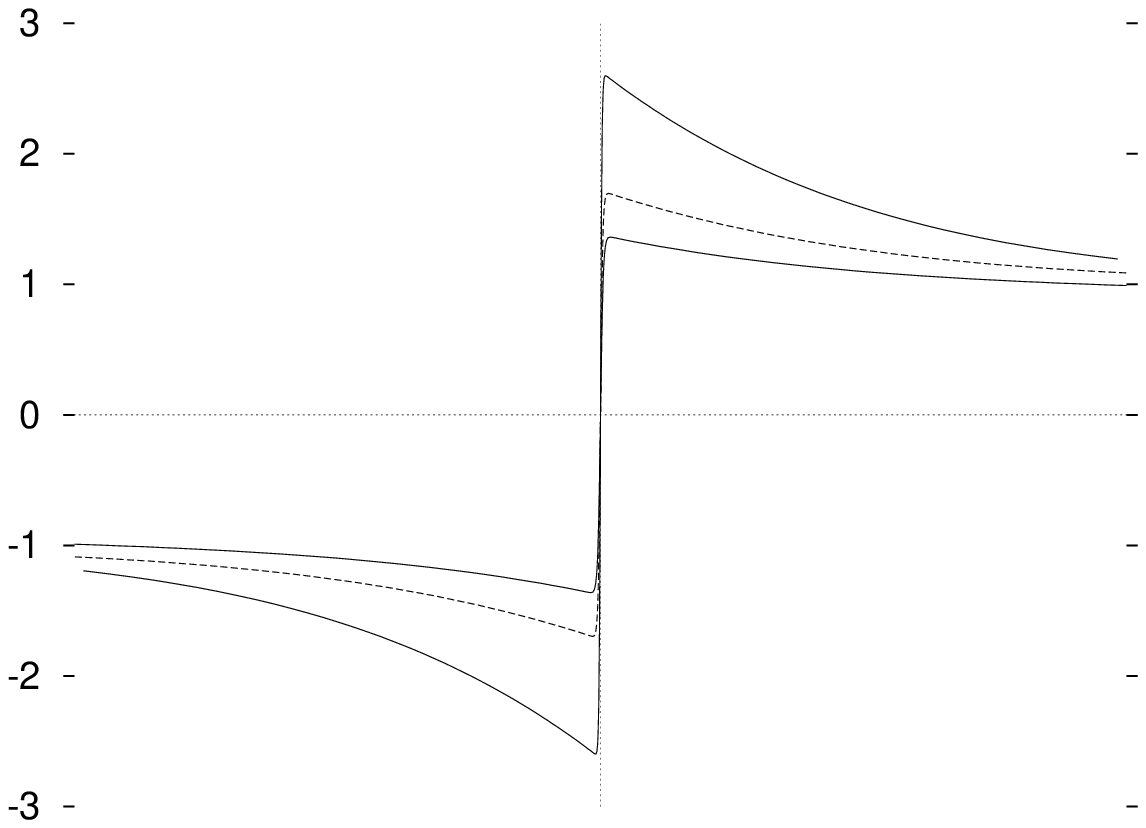}}
\subfigure[$V$-component]{\includegraphics[scale=.4]{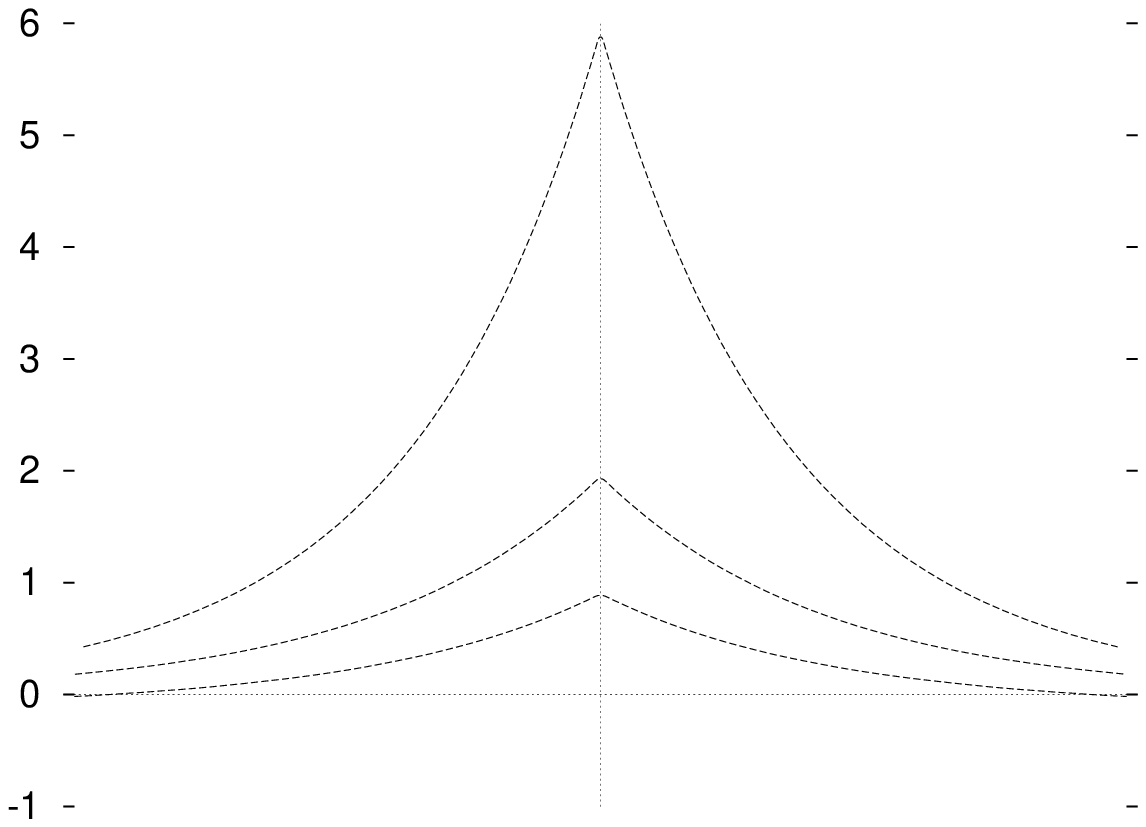}}
\caption{Numerical simulation of destabilization caused by the
  discrete spectrum, both components blow up in finite time ($H_0=-1$,
  $\ga=1.4$, $\tau=1$ and $\eps=0.1$).}
\end{center}
\end{figure}

\begin{figure}[htb]\label{essblowup}
\begin{center}
\psfrag{x}{$\xi$}
\psfrag{t}{$t$}
\psfrag{u}{$U$}
\psfrag{v}{$V$}
   \subfigure[$U$-component]{\includegraphics[scale=.4]{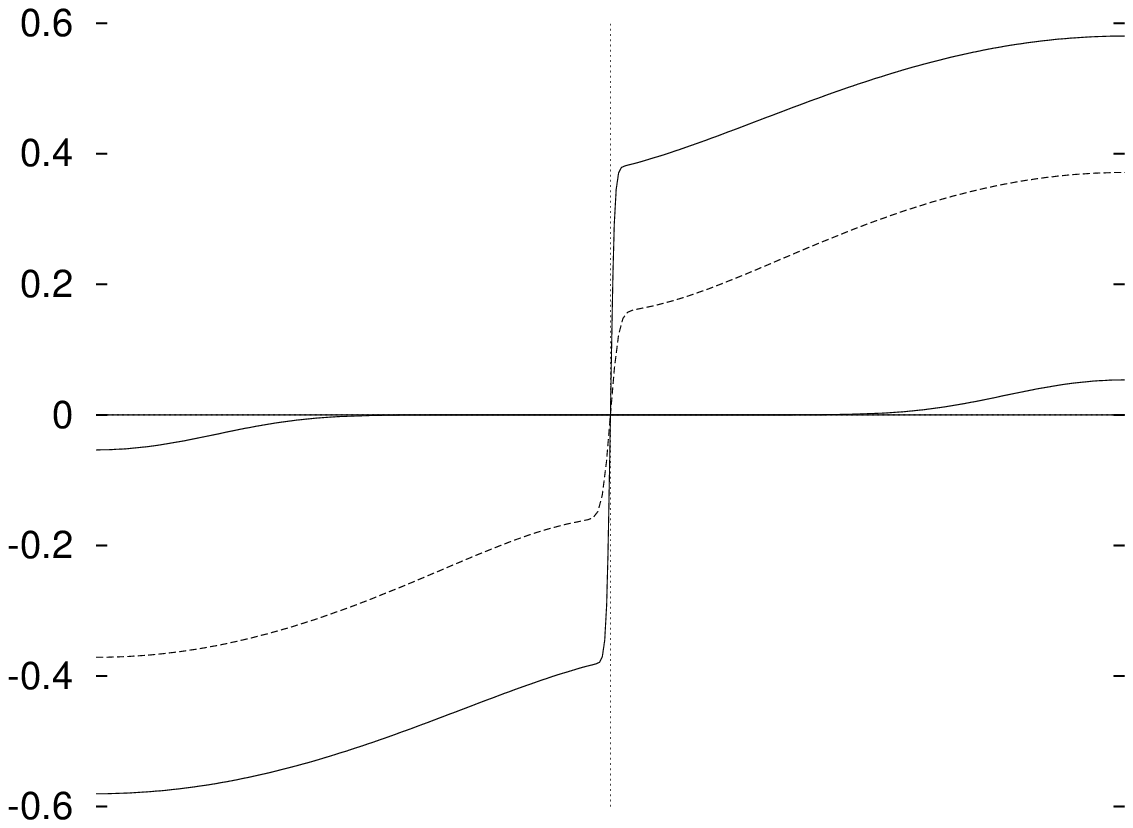}}
   \subfigure[$V$-component]{\includegraphics[scale=.4]{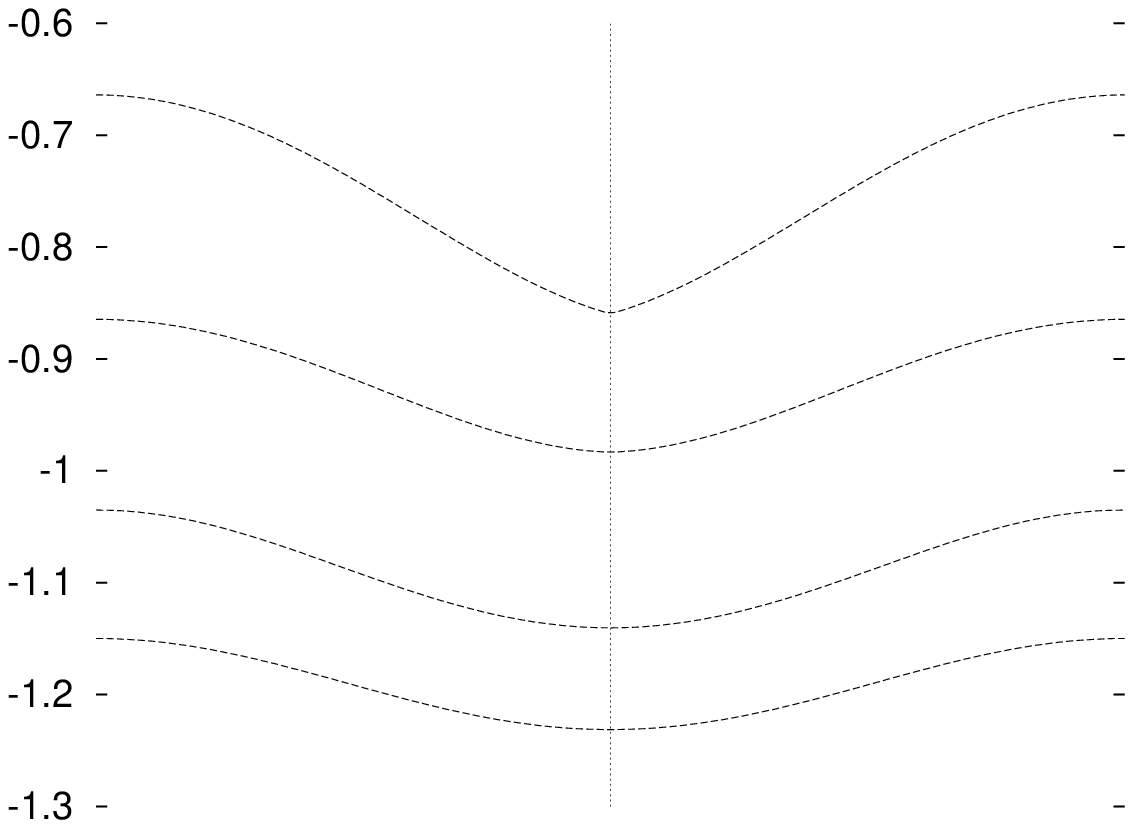}}
 \caption{Numerical simulation of destabilization caused by the
   essential spectrum, $U\to 0$ and $|V|$ grows slowly and
   exponentially ($H_0=1$, $\ga=-0.1$, $\tau=1$ and $\eps=0.1$).}
 \end{center}
 \end{figure}

We now examine numerically the difference between the two types of
bifurcations discussed in Theorems \ref{th:edge_SN_ex} and \ref{th:destab}. 
We consider the example system of sections \ref{sec:ex_suslo} and 
\ref{sec:stab_suslo_ex} for $H_0 > 0$ (case {\bf (i)}, type D) and 
$H_0 < 0$ (case {\bf (ii0}, type E). First, we note
that in both cases the simulations confirm that the fronts 
are asymptotically stable up to the 
analytically determined bifurcation values. In case {\bf (i)}
the front destabilizes at $\ga<\ga_{\rm double}$ due to an eigenvalue
in the discrete spectrum.  The eigenfunction associated to this type
of destabilization is localized to a neighborhood of the front as can
be seen in figure \ref{desblowup}.  In this case the front becomes
unstable and blows up in finite time, while the background states
remain stable.  In case {\bf (ii)}, 
the tip of the essential spectrum becomes positive and the background
states become unstable
as $\ga$ passes through $0$.  
As can be seen in figure \ref{essblowup}, this
destabilization causes the front to collapse.  The $U$ component then
tends to $0$ on the entire real line and the $V$ component 
grows according to  $V_t = V_{xx} +\eps^2 |\gamma| V$.
Thus, we may conclude that 
type D or type E orbits indeed exhibit significantly different behavior
at the destabilization. 
These simulations were performed using SPMDF \cite{blze94}, with
Neumann boundary conditions at $x=\pm50$.  The initial conditions used
in figure (\ref{desblowup}) are given by 
$U(x,0)=u_0(x,\eps; v_1)$ (\ref{u0p0} and $V(x,0) = v_1 e^{-\eps\sqrt{|\ga}|x|}$
(as described in Theorem \ref{th:ex_suslo}).

\subsection{Hopf bifurcations}
As we have seen in section \ref{sec:stab_reg}, in general there can be
(complex) eigenvalues near the endpoints $\la^{\pm}(0)$ of $\si_{\rm
  ess}$. Thus, if we keep $G_1 < 0$ fixed at an $\O(1)$ value and
increase $H_0$ such that $H_0 + G_1 - 2\tau$ approaches $0$, we
encounter a similar issue as was studied in the previous section: will
the front be destabilized by $\si_{\rm ess}$ at $H_0 = 2 \tau -G_1$ or
(just) before that, by an eigenvalue?  In this case, the bifurcation
is of Hopf type, and it is not associated to the existence problem.
This problem can in principle be analyzed by the methods developed
here, i.e. by determining $t_2(\la, \eps)$ through $\Delta_{\rm slow}
v_{\xi}$ and $\Delta_{\rm fast} v_{\xi}$.  We have already mentioned
the new features of the measuring the slow `jump' $\Delta_{\rm slow}
v_{\xi}$ in section \ref{sec:stab_suslo_bif}.  Moreover, since the
bifurcation does not occur near $\la = 0$, we do not have an explicit
formula for $u_{\rm in}(\xi)$, like (\ref{uin}), and it is thus not
immediately clear whether it is possible to determine $\Delta_{\rm
  fast} v_{\xi}$. Note that this latter issue is solvable with the
hypergeometric functions method developed in \cite{dgk1,dgk3}.
Nevertheless, we do not go deeper into this subject here.

\subsection{Planar fronts and stripes}
A next step in the study of (planar) stripes, as mentioned in the
Introduction, is the stability analysis of planar fronts, i.e. the
analysis of the stability of the fronts $(U_h(\xi),V_h(\xi))$ with
respect to two-dimensional perturbations (thus, $(U_h(\xi),V_h(\xi))$
represents a planar front that has a trivial structure in the
$y$-direction). The methods developed here can be used to study this
problem (as is also suggested by \cite{dp} in which a similar problem
has been studied in a mono-stable Gierer-Meinhardt context). It should
be noted here that there are several papers in the literature that
consider the question of the (non-)persistence of the stability of
one-dimensional fronts as two-dimensional planar fronts (see for
instance \cite{omk,tn,kap97,ns}). The analysis in \cite{tn,ns} of a
class of singularly perturbed bi-stable systems shows that the planar
fronts considered there cannot be stable, while it is shown that
planar fronts can be stable in a more regular context in \cite{kap97}.
Thus, this is a nontrivial issue. Preliminary analysis of the front
solutions considered in this paper indicates that these solutions
remain stable as planar fronts in the regular case (i.e. as long as
$G_1 < 0$ and $\O(1)$). The analysis of the planar fronts, and their
spatially periodic counterparts, the stripe patterns, is the subject
of work in progress.

{\bf Acknowledgements.}  D.I. would like to thank NSERC for their
support by way of a post doctoral fellowship.  A.D. and D.I
acknowledge support of the `Research Training Network (RTN):
Fronts-Singularities' (RTN contract: HPRN-CT-2002-00274).


\begin{thebibliography}{99}
  
\bibitem{agj} J. Alexander, R.A. Gardner, C.K.R.T. Jones [1990], A
  topological invariant arising in the stability analysis of
  travelling waves, {\it J. Reine Angew. Math.} {\bf 410}, 167--212.

\bibitem{blze94} J.G. Blom, P.A. Zegeling [1994], Algorithm 731: A
  Moving-Grid Interface for Systems of One-Dimensional Time-Dependent
  Partial Differential Equations, {\it ACM Transactions in
    Mathematical Software}, {\bf 20}, 194--214. 
  
\bibitem{dgk1} A.  Doelman, R. A. Gardner and T.J. Kaper [1998],
  Stability analysis of singular patterns in the 1-D Gray--Scott
  model: A matched asymptotics approach, {\it Physica D} {\bf 122},
  1--36.
  
\bibitem{dgk2} A. Doelman, R. A. Gardner, and T.J. Kaper [2002], A
  stability index analysis of 1-D patterns of the Gray--Scott model,
  {\it Memoirs of the AMS} {\bf 155} (737).
  
\bibitem{dgk3} A. Doelman, R. A. Gardner, and T.J. Kaper [2001], Large
  stable pulse solutions in reaction-diffusion equations, {\it Ind.
    Univ. Math. J.}, {\bf 50}(1), 443--507.
  
\bibitem{dkz} A. Doelman, T.J. Kaper, and P. Zegeling [1997], Pattern
  formation in the one-dimensional Gray-Scott model, {\it
    Nonlinearity} {\bf 10}, 523--563.
  
\bibitem{dp} A. Doelman and H. van der Ploeg [2001], Homoclinic stripe
  patterns, {\it SIAM J. Appl. Dyn. Syst.} {\bf 1}, 65--104.
  
\bibitem{henry} D. Henry [1981], {\it `Geometric Theory of Semilinear
    Parabolic Equations'}, Lecture Notes in Math. {\bf 840},
  Springer-Verlag.
  
\bibitem{fen79} N. Fenichel [1979], Geometrical singular perturbation
  theory for ordinary differential equations, {\it J. Diff. Eq.}  {\bf
    31}, 53--98.
  
\bibitem{gj} R.A. Gardner and C.K.R.T. Jones [1991], Stability of the
  travelling wave solutions of diffusive predator-prey systems, {\it
    Trans. AMS} {\bf 327}, 465--524.
  
\bibitem{jon95} C.K.R.T. Jones [1995], Geometric singular perturbation
  theory, in {\it Dynamical systems, Montecatibi Terme, 1994}, Lecture
  Notes in Mathematics {\bf 1609}, R. Johnson (ed.), Springer-Verlag.
  
\bibitem{kap97} T. Kapitula [1997], Multidimensional stability of
  planar travelling waves, {\it Trans. Amer. Math. Soc.} {\bf 349},
  257--269.

\bibitem{kap98} 
T. Kapitula [1998], 
The Evans function and generalized Melnikov integrals, 
{\it SIAM J. Math. Anal.} {\bf 30}, 273--297.

\bibitem{ks}
T. Kapitula and B. Sandstede [2002],
Edge bifurcations for near integrable systems via Evans function
techniques, {\it SIAM J. Math. Anal.} {\bf 33}, 1117--1143.

\bibitem{ns}
Y. Nishiura and H. Suzuki [1998],
Nonexistence of higher dimensional stable Turing patterns in the singular limit,
{\it SIAM J. Math. Anal.} {\bf 29}, 1087--1105.

\bibitem{omk}
T. Ohta, M. Mimura and R. Kobayashi [1989],
Higher-dimensional localized patterns in excitable media,
{\it Physica D} {\bf 34}, 115--144.

\bibitem{pw}
R.L. Pego, M.I. Weinstein [1992],
Eigenvalues, and instabilities of solitary waves,  
{\it Philos. Trans. Roy. Soc. London Ser. A} {\bf 340}, 47--94.

\bibitem{rob}
C. Robinson [1983],
Sustained resonance for a nonlinear system with slowly-varying
coefficients,
{\it SIAM J. Math. Anal.} {\bf 14}, 847--860.

\bibitem{tn}
M. Taniguchi and Y. Nishiura [1994],
Instability of planar interfaces in reaction-diffusion systems,
{\it SIAM J. Math. Anal.} {\bf 25}, 99--134.

\bibitem{titch}
E.C. Titchmarsh [1962],
{\it Eigenfunction Expansions Associated with Second-order Differential
Equations} (2nd ed.), Oxford Univ. Press.

\end{thebibliography}
\end{document}